\documentclass{amsart}

\usepackage{graphicx}

\def\R{{\mathbb R}}

\def\N{{\mathbb N}}

\def\Z{{\mathbb Z}}
\def\1{{1\!\!\!1}}
\def\a{{\alpha}}
\def\E{{\mathbb E}}
\def\P{{\mathbb P}}

\def\cal{\mathcal}

\def\ol{\overline}

\def\dist{{\rm{dist}}}

\def\supp{{\rm{supp}}}
\def\eps{\varepsilon}

\newcommand{\be}{\begin{equation}}
\newcommand{\ee}{\end{equation}}
\numberwithin{equation}{section}

\newtheorem{theorem}{Theorem}
\newtheorem{prop}{Proposition}[section]
\newtheorem{cor}{Corollary}[section]

\newtheorem{lemma}{Lemma}[section]

\title{Martin boundary of a reflected random walk on a half-space}
\author{Irina Ignatiouk-Robert}
\address{
{Universit\'e de Cergy-Pontoise,}
{D\'epartement de math\'ematiques,}
{2, Avenue Adolphe Chauvin,}
{95302 Cergy-Pontoise Cedex,}
{France}}
\date{\today}
\email{Irina.Ignatiouk@math.u-cergy.fr}
\keywords{Martin boundary.Sample path large deviations. Random walk.} 
\subjclass{Primary 60J50;  Secondary 60F10; 60J45}

\begin{document}
\begin{abstract} The complete representation of the 
Martin compactification for reflected random walks on a half-space $\Z^d\times\N$ is
obtained. It is shown that the full Martin compactification  is in
general not homeomorphic to the ``radial'' compactification obtained by Ney and Spitzer
for the homogeneous random walks in $\Z^d$ :   convergence of a sequence of points
$z_n\in\Z^{d-1}\times\N$ to a point of on the Martin boundary does not imply  convergence of
the sequence $z_n/|z_n|$ on the unit sphere $S^d$. Our approach relies on the  large
deviation properties of the 
scaled processes and uses Pascal's method combined with the ratio limit theorem. The
existence of non-radial limits is related to non-linear optimal large deviation
trajectories. 
\end{abstract}
\maketitle

\section{Introduction and main results}\label{sec1}

For an irreducible transient Markov chain  $(Z(t))$ on a countable set $E$ having Green's function
$G(z,z')$, the Martin compactification $\tilde{E}$ is
 the smallest compactification of the set $E$ for which the Martin kernels 
\[
K(z,z') = G(z,z')/G(z_0,z') 
\]
extend continuously with respect to the second variable $z'$ for every $z\in E$. A point
$\eta\in\partial E = \tilde{E}\setminus E$ is said to belong to the minimal Martin
boundary if  $K(\cdot, \eta)$ is a minimal harmonic function~(see Woess~\cite{Woess} for
the precise definitions). An explicit
representation of the Martin boundary and the minimal Martin boundary $\partial_M E
\subset \partial E$ allows to describe all harmonic
functions of the Markov chain $(Z(t))$~: by Poisson-Martin representation theorem, every
positive harmonic function $h$ 
is of the form
\[
h(z) = \int_{\partial_M E} K(z,\eta) \, d\nu(\eta) 
\] 
where $\nu$ is a positive Borel measure on $\partial_M E$. Moreover, by
convergence theorem, for every $z\in E$, the 
sequence $Z(n)$ converges 
$P_z$ almost surely to a random variable taking the values in $\partial_M E$.  

An explicit
description of the Martin boundary is usually a non-trivial problem. The most of the
existing results in this domain were obtained for the homogeneous processes (see
Woess~\cite{Woess} and the references
therein). One of the few results where the full Martin compactification was
obtained for non-homogeneous processes is the paper of Kurkova and
Malyshev~\cite{Kurkova-Malyshev}. They considered  random walks on a
half-plane $\Z\times\N$ and in the quadrant $\Z^2_+ = \N\times\N$ which behave as a
homogeneous nearest neighbors random walk in the interior of the domain and have some
different (homogeneous) transition probabilities on the boundary. Their results 
show a very surprising relationship between the Martin 
compactification and the optimal large deviation trajectories described  for such
processes obtained in Ignatyuk, Malyshev and 
Scherbakov~\cite{I-M-S}. Let us illustrate this
relationship on the example of the reflected random walks on the half-plane : the results of Kurkova
and Malyshev show that for such a random walk, there are two real values $0 \leq \theta_1
\leq\theta_2 \leq \pi$ such that 
\begin{itemize}
\item[i)] a sequence of points $z_n \in\Z\times\N$ with $\lim |z_n| = \infty$ converges
  to a point $\eta(\theta)$  of the Martin boundary if the sequence $z_n/|z_n|$ converges to a point
  $e^{i\theta}$ on $S^2_+ =\{ e^{i\theta} : \theta\in[0,\pi]\}$; 
\item[ii)] two sequences $z_n, z_n' \in\Z\times\N$ with $\lim |z_n| = \lim |z_n'| =
  \infty$, $\lim z_n/|z_n| = e ^{i\theta}$ and $\lim z_n'/|z_n'| = e ^{i\theta'}$ converge
  to the same point $\eta(\theta)=\eta(\theta')$ of the Martin boundary if an only if 
\begin{itemize}
\item[--] either $\theta = \theta' \in [\theta_1,\theta_2]$, \; $\mod(2\pi)$, 
\item[--] or $\theta, \theta' \in [0,\theta_1]$, \; $\mod(2\pi)$, 
\item[--] or $\theta, \theta'\in [\theta_2,\pi]$, \; $\mod(2\pi)$. 
\end{itemize} 
\end{itemize}
In ~\cite{I-M-S} it was shown that for every $T>0$, the family of scaled random walks
$$(Z^\eps(t) ~\dot=~ \eps Z([t/\eps]), \, t\in[0,T])$$  satisfy sample path large deviation
principle with a rate function $I_{[0,T]}(\phi)$ and that  with the same values $\theta_1$
and $\theta_2$, the following assertions hold.  
\begin{itemize}
\item[--] For $\theta\in[\theta_1,\theta_2]$,  the optimal large deviation trajectory 
$\phi_\theta :[0,T_\theta]\to\R\times\R_+$ minimizing the  rate function
$I_{[0,T]}(\phi)$ over all $T>0$ and all continuous functions $\phi :
[0,T]\to\R\times\R_+$ with given $\phi(0)=0$ and $\phi(T)= e^{i\theta}\in S^2_+$  is linear :
$\phi_\theta(t) = e^{i\theta} t/T_\theta$ with some $T_\theta >0$; 
\item[--] while for $\theta\in[0,\theta_1]\cup[\theta_2,\pi]$, such a trajectory is piece-wise linear and is of
  the form  
\[
\phi_\theta(t) ~=~ \begin{cases} \gamma_\theta t/T_\theta' &\text{ for $t\in[0,T_\theta']$}\\  
\gamma_\theta + (e^{i\theta} - \gamma_\theta) (t-T_\theta')/(T_\theta - T_\theta')&\text{
  for $t\in[T_\theta',T_\theta]$}   
\end{cases}
\]
with some $T_\theta > T_\theta > 0$ where 
$\gamma_\theta$ is a unique point on the boundary $\R\times\{0\}$
for which 
\[
\arg( e^{i\theta} - \gamma_\theta) = \begin{cases} \theta_1 &\text{ if $\theta\in[0,\theta_1]$,}\\ 
\theta_2 &\text{ if $\theta\in[\theta_2,\pi]$.}
\end{cases}
\]
\end{itemize}
Unfortunately, the method proposed by Kurkova and Malyshev~\cite{Kurkova-Malyshev}
required very particular properties of 
the process : they considered the random walks for
which the  only non-zero transitions in  the interior of 
the domain are on the nearest neighbors:  $p(z,z \pm e_i) = \mu(\pm e_i)$ with $e_1=(1,0)$
and  $e_2=(0,1)$. For such  random walks,  the jump  generating function  is defined  by $$
\varphi(x,y) ~=~  \mu(e_1) x  + \mu(-e_1)x^{-1}  + \mu(e_2)y +  \mu(-e_2)y^{-1} $$  and the
equation $  xy(1-\varphi(x,y)) ~=~  0 $ determines  an elliptic  curve ${\bf S}$  which is
homeomorphic to  the torus.   To identify  the Martin boundary,  a functional  equation
was 
derived for  the generating function  of the Green's  function and the asymptotics  of the
Green's function were calculated by using the methods of complex  analysis on the elliptic
curve ${\bf S}$.  Such a method seems  to be unlikely to apply in a more general situation,
for  higher dimensions  or when  the  jump sizes  are arbitrary,  because the proof is
based on the geometrical
properties of the elliptic curve ${\bf S}$ : even for the $2$-dimensional case, if a
random walk has an additional non-zero 
transition $p(z,z + u) = \mu(u)$ with $|u| > 2$,  the equation $  xy(1-\varphi(x,y)) ~=~
0 $ is not of the second order and consequently, the corresponding elliptic curve is not
homeomorphic to the torus. 

Since the large deviation methods extend easily for an arbitrary dimension and for
arbitrary jumps, a natural idea is to use them in order to identify the Martin
boundary. The similarities of the results of Kurkova
and Malyshev~\cite{Kurkova-Malyshev} and the large deviation results of  Ignatyuk, Malyshev and
Scherbakov~\cite{I-M-S} suggest that such an approach should be possible. The first result
in this domain was
obtained in Ignatiouk-Robert~\cite{Ignatiouk:06} for a homogeneous random walk $(Z_+(t))$ on $\Z^d$
killed upon hitting the negative half-space $\Z^{d-1}\times(-\N)$~: the large deviation
technique was combined there with Bernoulli part decomposition due to Foley and
McDonald~\cite{Foley-McDonald}. The main steps of this method can be summarized as
follows~: 
\begin{itemize}
\item[--] The first step is a  ratio limit theorem:~ Bernoulli part decomposition was used
  to identify the limits of the Martin kernel 
$K(z,z_n)$ when the logarithmic asymptotic of Green's function for a given sequence $(z_n)$ is
zero. 
\item[--] The logarithmic asymptotics of Green's function were obtained with the large
  deviation technique. 
\item[--] An appropriated exponential change of the measure was finally used in order to apply the
  ratio limit theorem for a twisted Markov process for which the corresponding logarithmic asymptotic of Green's
  function is zero.  
\end{itemize}
In the present paper the large deviation method is developed in order to identify the
Martin boundary for a reflected random walk $(Z(t))$ on the 
half-space $\Z^{d-1}\times\N$. Such a random walk behaves as a homogeneous random walk in
the interior of the half-space and has some different transition probabilities on the
boundary hyper-plane $\Z^{d-1}\times\{0\}$. Here, the approach of 
Ignatiouk-Robert~\cite{Ignatiouk:06} is not only harder to apply but 
also it does not work in general because the corresponding twisted process does not
exist. To solve this problem we refine the large deviation technique.  

 We show that the family of scaled processes $(Z^\eps(t)=\eps
  Z([t/\eps]), t\in[0,T])$ satisfies sample path large deviation principle with a good rate function
  $I_{[0,T]}$ and that the logarithmic asymptotics of  Green's function $G(z,z_n)$
   of the original process $(Z(t))$ when $|z_n|\to\infty$ and $z_n/|z_n|\to q$ are
   determined by the quasi-potential 
\[
I(0,q) ~=~ \inf_{T>0}~\inf_{\phi :~\phi(0)=0, \; \phi(T)=q} I_{[0,T]}(\phi) 
\]
which represents an optimal large deviation
cost to go from the point $0$ to the point $q$. Next, the method of
\cite{Ignatiouk:06} is used to identify the limit of the Martin kernel 
$K(z,z_n)$ when $|z_n|\to\infty$ and the limit $z_n/|z_n|\to q$ belongs to the boundary hyper-plane
$\R^{d-1}\times\{0\}$. This is the first step of our proof.

For $q\not\in\R^{d-1}\times\{0\}$ we consider a function $\phi:
  [0,T]\to\R^d$ with $\phi(0)=0$ and $\phi(T)=q$ where the minimum $I(0,q)$ is achieved. Such a function $\phi$ represents an 
  optimal large deviation path from $0$ to $q$. It is shown that every optimal  large
  deviation path from $0$ to $q$ leaves the boundary hyper-plane $\R^{d-1}\times\{0\}$ at
  some point $\gamma_q\in \R^{d-1}\times\{0\}$ and that Green's function $G(z,z_n)$ can be
  decomposed into a main part determined by $\gamma_q$ 
  and the corresponding negligible part. The main part of $G(z,z_n)$ corresponds to  the trajectories of the process $(Z(t))$
  that leave the boundary hyper-plane in a $\delta|z_n|$-neighborhood of the point
  $\gamma_q|z_n|$. With this approach we identify the  limit of the Martin kernel 
$K(z,z_n)$ when $|z_n|\to\infty$ and  $z_n/|z_n|\to q$ for any
  $q\in\R^{d-1}\times[0,+\infty[$. 

The reflection on the boundary is not only harder to tackle but also yields very different
and interesting 
results. Contrary to the case analyzed in
\cite{Ignatiouk:06}, here the convergence to the Martin boundary can be
non-radial~: a convergence to a point on the Martin boundary of a sequence $(z_n)$ does
not imply the convergence of the sequence $z_n/|z_n|$ on the unit sphere. We obtain this
result as a
consequence of the existence of non-linear optimal large deviation trajectories. 
\subsection{Main result} 
We consider a Markov process $Z(t)=(X(t),Y(t))$ on
$\Z^{d-1}\times\N$ with transition probabilities 
\be\label{e1-1}
p(z,z') ~=~ \begin{cases} \mu(z'-z) &\text{for ~$z=(x,y),z'\in\Z^{d-1}\times\N$ with $y>0$,}\\
\mu_0(z'-z) &\text{for $z=(x,y),z'\in\Z^{d-1}\times\N$ with $y=0$} 
\end{cases}
\ee
where $\mu$ and $\mu_0$ are two different probability measures on $\Z^d$ having the means 
\be\label{e1-2}
m \dot= \sum_{z\in\Z^d} z \mu(z)  \quad \quad \text{and} \quad \quad m_0 \dot= \sum_{z\in\Z^d} z
\mu_0(z). 
\ee
Throughout this
paper we denote by $\N$ the set of all non-negative integers~: $\N=\{0,1,2,\ldots\}$ and
we let $\N^*=\N\setminus\{0\}$. The assumptions we need on the Markov process
$(Z(t))$ are the following.

\medskip
\noindent 
\begin{enumerate}
\item[(H0)] {\em $\mu(z)=0$ for $z=(x,y)\in\Z^{d-1}\times\Z$ with $y<-1$ and $\mu_0(z) = 0$ for
$z=(x,y)\in\Z^{d-1}\times\Z$ with $y< 0$.} 
\item[(H1)] {\em The Markov process $Z(t)$ is irreducible on
  $\Z^{d-1}\times\N$.} 
\item[(H2)] {\em The homogeneous random walk $S(t)$ on $\Z^d$ 
having transition probabilities  
$p_S(z,z')=\mu(z'-z)$ is irreducible on $\Z^d$ and the last coordinate of $S(t)$ is an aperiodic random walk on $\Z$ .}
\item[(H3)]  
\[
m \not= 0 \quad \quad \text{ and } \quad \quad \frac{m}{|m|} + \frac{m_0}{|m_0|} \not= 0.
\]
\item[(H4)]{\em  The jump generating functions 
\be\label{e1-3}
\varphi(a) ~=~ \sum_{z\in\Z^d} \mu(z) e^{a\cdot z} \quad \text{ and } \quad \varphi_0(a) ~=~
\sum_{z\in\Z^d} \mu_0(z) e^{a\cdot z} 
\ee
are finite everywhere on $\R^d$.} 
\end{enumerate}
Under the above assumptions, the sets 
\be\label{e1-4}
D ~\dot=~ \{a\in \R^d :
\varphi(a) \leq 1\} \quad \text{and} \quad D_0 ~\dot=~ \{a\in \R^d :
\varphi_0(a) \leq 1\} 
\ee 
are convex and the set $D$ is moreover compact (see~\cite{Hennequin}). 
The following parts
of  the boundary $\partial D$
are important for our analysis :  
\[
\partial_0 D ~\dot=~ \{a\in\partial D :~ \nabla\varphi(a) \in
\R^{d-1}\times\{0\}\}
\] 
\[
\partial_+ D ~\dot=~ \{a\in\partial D :~ \nabla\varphi(a) \in
\R^{d-1}\times[0,+\infty[\}
\] 
and  
\[
\partial_- D ~\dot=~ \{a\in\partial D :~ \nabla\varphi(a) \in
\R^{d-1}\times ]-\infty,0]\}. 
\]
 For   
$a\in D$,  denote by $\ol{a}$  the unique
point on the boundary $\partial_- D$  which has the same first $(d-1)$ coordinates as the
point $a$ and let
\be\label{e1-5}
\hat{D} = \{a\in D : \varphi_0(\ol{a})\leq 1\}. 
\ee
Remark that under the hypotheses~(H0)-(H1), for any $a\in D$,
\[
\varphi_0(\ol{a}) \leq \varphi_0(a) 
\]
because the function $a=(\a,\beta)\to\varphi_0(a)$ is increasing with respect to the last coordinate
$\beta$ of $a=(\a,\beta)\in\R^d$. This inequality implies  another useful
representation of the set $\hat{D}$ : 
\[\text{ $a=(\a,\beta)\in\hat{D}$ if
and only if $a\in D$ and $a'=(\a,\beta')\in D\cap D_0$ for some  $\beta' \in\R$ } 
\]
or equivalently,
\be\label{e1-6}
\hat{D} ~=~ (\Theta\times\R) \cap D 
\ee
where 
\be\label{e1-7}
\Theta ~\dot=~ \{ \a \in\R^{d-1} :~ \inf_{\beta\in\R} \max\{\varphi(\a,\beta),
\varphi_0(\a,\beta)\} \leq 1\}.
\ee 
The set $\Theta\times\{0\}$ is therefore the orthogonal projection of the set $D\cap D_0$ onto the hyper-plane
$\R^{d-1}\times\{0\}$. Remark finally that $\partial_0 D = \partial_+ D\cap \partial_-D$ and for $a\in\partial_+
D$,  $a=\ol{a}$
if and only if $a\in\partial_0 D$. 

It is moreover convenient to introduce  the following notations~:~ for 
$a\in \hat{D}=(\Theta\times\R) \cap D$, we denote by $V(a)$ the normal cone to the set 
$\hat{D}$ at the point $a$ and for $a\in\hat{D}\cap\partial_+D= (\Theta\times\R) \cap \partial_+D$ we define the function
$h_a$ on $\Z^{d-1}\times\N$ by letting 
\be\label{e1-8}
h_a(z) ~=~ \begin{cases} \displaystyle{e^{a\cdot z}  -
  ~\frac{1- \varphi_0(a)}{1-
  \varphi_0(\ol{a})} \, e^{\ol{a}\cdot z}}  &\text{ if  
 $\; a\not\in\partial_0 D \; $ and  $\; \varphi_0(\ol{a}) < 1$,}\\ 
\\
\displaystyle{y e^{a\cdot z} +  
 ~\frac{ \frac{\partial}{\partial\beta}\varphi_0(a)}{(1-
  \varphi_0(a))} \, e^{a\cdot z}} &\text{ if  $\; a=\ol{a}\in\partial_0 D \; $ and $\; \varphi_0(a) < 1$,}\\
\\
e^{\ol{a}\cdot z} &\text{ if $\; \varphi_0(\ol{a}) = 1\; $}
\end{cases}
\ee
 where $\frac{\partial}{\partial\beta}\varphi(a)$ denotes the partial derivative of the
function $a\to \varphi(a)$ with respect to the last coordinate $\beta\in\R$ of $a=(\a,\beta)$.   

\medskip 

We denote by ${\cal S}_+^d$ a half-sphere ${\cal S}^d\cap\R^{d-1}\times\R_+$ and
$G(z,z')$ denotes Green's function of the Markov process $(Z(t))$. 

\medskip

Our preliminary results show that for
any $q\in{\cal S}^d_+$, there is a unique point
  $\hat{a}(q)\in \hat{D}\cap\partial_+D$ for which $q\in V(\hat{a}(q))$ 
and that for every $a\in \hat{D}\cap\partial_+D$, 
\be\label{e1-9}
V(a) ~=~ \begin{cases} \bigl\{c\nabla\varphi(a) : c\geq 0\bigr\} &\text{if
    either 
    $\varphi_0(\ol{a}) < 1$} \\
&\text{or $a=\ol{a}\in\partial_0 D$,}\\
\bigl\{c_1\nabla\varphi(a) + c_2(\nabla\varphi_0(\ol{a}) + \kappa_{a}
\nabla\varphi(\ol{a})) :~ c_i\geq 0\bigr\} &\text{if
    $\varphi_0(\ol{a}) = 1$ }\\ &\text{and $a\not\in\partial_0 D$}\\ 
\end{cases}
\ee
where 
\[
\kappa_{a} = - \left.\frac{\partial\varphi_0(\a,\beta)}{\partial\beta}
\left(\frac{\partial\varphi(\a,\beta)}{\partial\beta}\right)^{-1}\right|_{(\a,\beta) = \ol{a}} 
\]
(see Lemma~\ref{lem2-3}  and Lemma~\ref{lem2-5} below). 

\medskip

The main result of our paper is the following theorem. 

\begin{theorem}\label{th1-1} Under the hypotheses (H0)-(H4), the following assertions hold~: 
\begin{itemize}
\item[(i)] the Markov process $Z(t)$ is
  transient; 
\item[(ii)] for any
  $a\in\hat{D}\cap\partial_+ D$ and any sequence of points
  $z_n\in\Z^{d-1}\times\N$ with $\lim_{n\to\infty} |z_n|=\infty$, 
\be\label{e1-10}
\lim_{n\to\infty} G(z,z_n)/G(z_0,z_n) = h_a(z)/h_a(z_0),\quad \quad \quad \forall \;
z\in\Z^{d-1}\times\N 
\ee 
when $\lim_{n\to\infty}
  \dist(V(a),z_n/|z_n|)=0$.
\end{itemize}
\end{theorem}
\noindent
Assertion (ii) proves  that a sequence  $z_n\in\Z^{d-1}\times\N$ with
$\lim_{n\to\infty}|z_n| = \infty$, converges to a point on the Martin boundary  if and only if 
\[
\lim_{n\to\infty} ~\dist\left(V(a),z_n/|z_n|\right)=
  0 
\]
for some $a\in\hat{D}\cap\partial_+ D$. 
Recall that for a homogeneous random walk on $\Z^d$ (see Ney and
Spitzer~\cite{Ney-Spitzer}), a sequence $z_n\in\Z^d$ converges to
a point of the Martin boundary if and only if $\lim_{n\to\infty}|z_n| = \infty$ and the
sequence $z_n/|z_n|$ converges to a point on the unit sphere $S^d$. For the reflected random
walk on the half-space $\Z^{d-1}\times\N$, Theorem~\ref{th1-1} provides the
  existence of   non-radial limits : if the mapping $\hat{a} : {\cal S}_+^d\to
  \hat{D}\cap\partial_+ D$ is not one to one then the convergence to a point on the
  Martin boundary does not imply convergence of the sequence $z_n/|z_n|$. The explicit
  representation \eqref{e1-9} of the normal cone $V(a)$ shows that such a
  mapping is not one to one in a quite general situation~: when $\varphi_0(\ol{a}) = 1$ for
  some $\ol{a}\in\partial_-D$.

\subsection{The overview of the proof} To prove Theorem~\ref{th1-1} we identify first the
  harmonic functions of the process $(Z(t))$. Since the transition probabilities of the
  Markov process $(Z(t))$ are invariant with respect to the translations on
  $z\in\Z^{d-1}\times\{0\}$ and since the Markov process $(Z(t))$ is irreducible then the
  same arguments as in  Doob, Snell and Williamson~\cite{Doob:03} (see the proof of
  Theorem~5) show that every minimal harmonic function is of the form 
\[
h(x,y) = \exp(\a\cdot x) h(0,y), \quad \quad \forall (x,y)\in\Z^{d-1}\times\N
\]
with some $\a\in\R^{d-1}$. We prove that the constant multiples of the functions $h_a$
with  $a=(\a,\beta)\in\hat{D}\cap 
  \partial_+ D$, are the only minimal non-negative harmonic
  functions of the Markov process $(Z(t))$.  These
  arguments prove the first assertion of Theorem~\ref{th1-1} because under our hypotheses,
  $\{0\} \subset \hat{D}\cap\partial_+ D\not=\{0\}$.

To prove the assertion (ii), we identify first the logarithmic asymptotics of Green's
function, by using the large deviation method.   
The results of Dupuis, Ellis and Weiss~\cite{D-E-W}, Dupuis and
Ellis~\cite{D-E} and Ignatiouk~\cite{Ignatiouk:02,Ignatiouk:04} are used to show that  the family of  scaled 
  processes $(Z^\eps(t) = \eps Z([t/\eps]), \, t\in[0,T])$ satisfies sample path large deviation principle with a
  good rate function  $I_{[0,T]}(\phi)$ having an explicit form. The quasi-potential
  $I(0,q)$ of the rate function $I_{[0,T]}(\phi)$ represents an optimal large deviation
  cost to go from the point $0$ to the point $q$~: 
\[
I(0,q) ~=~ \inf_{T>0} ~\inf_{\phi : \phi(0)=0, \, \phi(T)=q} I_{[0,T]}(\phi)
\]
We show that  for any
$q\in\R^{d-1}\times\R_+$ and any sequence of points $z_n\in\Z^{d-1}\times\N$ with $\lim
|z_n|=+\infty$ and $\lim z_n/|z_n| =q$,
the following equalities hold 
\be\label{e1-11}
\lim_{n\to \infty} ~\frac{1}{|z_n|}\log G(z,z_n) ~=~ - I(0,q) ~=~ \sup_{a\in
  \hat{D}} ~a\cdot q ~=~ \hat{a}(q)\cdot q, \quad \quad \forall
z\in\Z^{d-1}\times\N.
\ee
When $\lim_{n\to\infty} z_n/|z_n| =q\in\R^{d-1}\times\{0\}$ and $\ol{\hat{a}(q)} = 0$, 
the proof of \eqref{e1-10} uses the following arguments : 
\begin{itemize}
\item[--] from \eqref{e1-11} we obtain the equality  
\[
\lim_{n\to\infty}~\frac{1}{n} \log G\bigl(z_0, z_n\bigr)   ~=~ 0 
\]
\item[--] and next, using the ratio limit
theorem of~\cite{Ignatiouk:06} we get \eqref{e1-10}.
\end{itemize} 

To get  \eqref{e1-10} for a sequence $z_n\in\Z^{d-1}\times\N$ with $\lim_{n\to\infty} z_n/|z_n|
=q\in\R^{d-1}\times\{0\}$ and $\ol{\hat{a}(q)} \not= 0$, the above arguments are combined
together with the exponential change of measure : the ratio limit theorem is applied for
a sub-stochastic twisted Markov chain having transition probabilities $
\tilde{p}(z,z') ~=~ p(z,z') \exp(a \cdot(z'-z))$
with a parameter $a=\ol{\hat{a}(q)}$. 

Similar arguments are used  in order to prove \eqref{e1-10} for a sequence of points
$z_n\in\Z^{d-1}\times\N$ with $\lim_{n\to\infty} z_n/|z_n| 
=q\in\R^{d-1}\times\R_+^*$ when 
$\varphi_0\left(\ol{\hat{a}(q)}\right) < 1$. 
The only difference  is here that  there is no suitable exponential change of
measure. Instead of the exponential change of measure we  
consider a twisted Markov chain  $(\tilde{Z}(t))$ with transition probabilities 
$
\tilde{p}(z,z') ~=~ p(z,z') h_{\hat{a}(q)}(z')/h_{\hat{a}(q)}(z)$. 
For the twisted Green's function $\tilde{G}(z,z') = G(z,z')
h_{\hat{a}(q)}(z')/h_{\hat{a}(q)}(z)$, 
the equality 
\[
\lim_{n\to\infty}~\frac{1}{n} \log \tilde{G}\bigl(z_0, z_n\bigr)   ~=~ 0
\]
follows from the relations \eqref{e1-11}  and the explicit
form of the harmonic function $h_{\hat{a}(q)}$. 

The  case  when $\lim_{n\to\infty} z_n/n = q\in\R^{d-1}\times\R_+^*$ and  
$\varphi_0\left(\ol{\hat{a}(q)}\right) = 1$ is more difficult.  
In this case, we can not use the above arguments  because there is no harmonic functions
$h$ satisfying the equality 
\[
\lim_{n\to\infty}~\frac{1}{n} \log \bigl(G\bigl(z_0, z_n\bigr) h(z_n)\bigr)   ~=~
0.
\]
Instead,  we use Pascal's method combined with the
renewal equation 
\be\label{e1-12}
G(z,z_n) ~=~ G_+(z,z_n) ~+ \sum_{\substack{w\in E\setminus \Z^{d-1}\times\{0\}, \\w'\in
  \Z^{d-1}\times\N^*}} G(z,w) p(w,w') G_+(w',z_n).  
\ee 
$G_+(z,z')$ denotes here the mean number of visits of the point $z'$ starting from $z$
before hitting the boundary hyperplane $\Z^{d-1}\times\{0\}$. Here, the main ideas of our proof
are the following :

For every point
$q\in\R^{d-1}\times\R_+^*$, the normal cone $V(\hat{a}(q))$ is generated by the vectors
$\gamma_q$ and $q-\gamma_q$ with some uniquely defined
$\gamma_q\in\R^{d-1}\times\{0\}$.  In 
the large deviation scaling, the point $\gamma_q$  corresponds to an
optimal way  from $0$ to $q$. Because of the influence of the boundary, the optimal ways
are not linear, an optimal way from $0$ to $q\in\R^{d-1}\times\R_+^*$ follows first a
linear trajectory on the
boundary hyper-plane $\R^{d-1}\times\{0\}$ before hitting the point
$\gamma_q\in\R^{d-1}\times\{0\}$ and next follows another linear trajectory  from
$\gamma_q$ to $q$ in the interior of the half-space.  

The right hand side of the renewal equation 
\eqref{e1-12} is decomposed into a principal part 
\[
\Xi_\delta^q(z,z_n) ~=~ G_+(z,z_n)\1_{\{\gamma_q=0\}}  ~+ \hspace{-0.5cm}\sum_{\substack{w\in \Z^{d-1}\times\{0\} : |w -  \gamma_q|z_n|| < \delta |z_n|, \\ w'\in
  \Z^{d-1}\times\N^*}} G(z,w) p(w,w') G_+(w',z_n) 
\]
corresponding to the optimal large deviation way to go from $0$ to $q$, and the
 negligible part
\[
G_+(z,z_n)\1_{\{\gamma_q\not=0\}}  ~+ \hspace{-0.5cm}\sum_{\substack{w\in \Z^{d-1}\times\{0\} : |w -  \gamma_q|z_n|| \geq \delta |z_n|, \\ w'\in
  \Z^{d-1}\times\N^*}} G(z,w) p(w,w') G_+(w',z_n). 
\]
Next, for those
$q\in\R^{d-1}\times\R^*$ for which 
$\gamma_q=0$, we obtain \eqref{e1-10}  by using the results of ~\cite{Ignatiouk:06}. 
When $\gamma_q\not=0$,  the equality \eqref{e1-10} is obtained from the convergence  
\[
 G(z,w)/G(z_0,w) ~\to~
h_{\hat{a}(\gamma_q)}(z)/h_{\hat{a}(\gamma_q)}(z_0) 
\]
as $|w|\to\infty$ and $w/|w|\to\gamma_q$. 
We use here the fact that $\hat{a}(\gamma_q) = \hat{a}(q)$ and that $\gamma_q\in\R^{d-1}\times\{0\}$
(recall that for $q\in\R^{d-1}\times\{0\}$, the equality \eqref{e1-10} is 
proved by using the ratio limit theorem). 

Our paper is organized as follows. Section~\ref{sec-pr} is devoted to the preliminary
results. 
The harmonic functions of the
Markov process $(Z(t))$  are identified in
Section~\ref{sec-hf}. In Section~\ref{sec-cc} we prove that our Markov process satisfies
strong communication condition. This property is needed to establish sample path large
deviation principle for the family of scaled processes and also to apply the ratio limit
theorem.   Section~\ref{sec-ld} is devoted to large deviation results. In
Section~\ref{sec-reneq} we apply large deviation results to decompose the right hand side
of the renewal equation \eqref{e1-12} into a principal part and a negligible
part. Section~\ref{sec-rlth} is devoted to the ratio limit theorem. The proof
Theorem~\ref{th1-1} is given in Section~\ref{sec-proof}.

\section{Preliminary results}\label{sec-pr}
Let $\tau =\inf\{t\geq 1~:~ Z(t)\in\Z^{d-1}\times\{0\}\}$
denote the first time when the process $Z(t)$ returns to the boundary hyper-plane
$\Z^{d-1}\times\{0\}$. Recall that for $z,z'\in\Z^{d-1}\times\N^*$, 
\[
G_+(z,z') ~=~ \sum_{t=0}^\infty \P_z(Z(t) = z', \, \tau > t) 
\]
is Green's function of a homogeneous random walk $Z_+(t)$ on $\Z^{d-1}\times\N^*$ having transition probabilities
$p_S(z,z') = \mu(z'-z)$ and killed upon hitting the
half-space $\Z^{d-1}\times(-\N)$. The homogeneous random walk on $\Z^d$ having transition probabilities
$p_S(z,z') = \mu(z'-z)$, $z,z'\in\Z^d$, and its Green's function are denoted by $S(t)$ and
$G_S(z,z')$ respectively. On several occasions we will need the following relations. 
\begin{lemma}\label{lem2-1} Under the hypotheses (H1) and (H2), for any $a\in D$   
\[
G_+(z,z') ~\leq~ \exp(a\cdot(z-z')) G_S(0,0) \quad \quad \forall z,z'\in\Z^{d-1}\times\N^*.
\]
If moreover $\varphi_0(a)\leq 1$ then also 
\[
G(z,z') ~\leq~ \exp(a\cdot(z-z')) G(z',z') \quad \quad \forall z,z'\in\Z^{d-1}\times\N.
\]
\end{lemma} 
\begin{proof} Indeed, for $a\in D$, the exponential function $z\to\exp(a\cdot z)$ is
  super-harmonic for the Markov process $Z_+(t)$. By Harnack's inequality from this it
  follows that 
\begin{align*}
G_+(z,z')/G_S(0,0) &~=~ G_+(z,z')/G_S(z',z') \\&~\leq~ G_+(z,z')/G_+(z',z') ~=~ \P_z(Z_+(t) = z'
\;\text{for some} \; t \in\N) \\&~\leq~ \exp(a\cdot(z-z'))  
\end{align*}
for all $z,z'\in\Z^{d-1}\times\N^*$. 
Moreover, for those $a\in D$ for which $\varphi_0(a)\leq 1$, the exponential function
$z\to\exp(a\cdot z)$  is also super-harmonic for the Markov process $Z(t)$. Hence, using
again Harnack's inequality we obtain  
\[
G(z,z')/G(z',z') ~=~ \P_z(Z(t) = z'
\;\text{for some} \; t \in\N) ~\leq~ \exp(a\cdot(z-z')) 
\]
for all $z,z'\in\Z^{d-1}\times\N$. 
\end{proof}

\begin{lemma}\label{lem2-2}  Under the hypotheses (H0) and (H1), for any $a\in D$ 
\[
\E_z\bigl(\exp(a\cdot Z(\tau)); \, \tau < \infty \bigr) ~=~ \begin{cases} \exp(\ol{a}\cdot z) &\text{if
      $z\in\Z^{d-1}\times\N^*$,}\\
\exp(\ol{a}\cdot z) \varphi_0(\ol{a}) &\text{if
      $z\in\Z^{d-1}\times\{0\}$}.
\end{cases}
\]
\end{lemma}
\begin{proof} Indeed, since $Z(\tau)\in\Z^{d-1}\times\{0\}$ then for any $a\in D$ and $z\in\Z^{d-1}\times\N$ 
\[
\E_z\bigl(\exp(a\cdot Z(\tau) - \ol{a}\cdot z); \, \tau < \infty\bigr) ~=~  \E_z\bigl(\exp(\ol{a}\cdot (Z(\tau)-z)); \,
\tau < \infty\bigr) 
\]  
because  according to the definition of the mapping
$a\to\ol{a}\in\partial_- D$ (see 
Section~\ref{sec1}), $a\cdot z = \ol{a}\cdot z$ for all $ z\in\Z^{d-1}\times\{0\}$.  
For $z\in\Z^{d-1}\times\N^*$, the right hand side of this equality 
is equal to the probability that a twisted homogeneous random walk $\tilde{Z}(t)$ on $\Z^d$ with
transition probabilities $\tilde{p}(z,z') ~=~ \exp(\ol{a}\cdot(z'-z))\mu(z'-z)$
starting from $z$ ever hits the  hyper-plane $\Z^{d-1}\times\{0\}$. Such a
twisted random walk has a finite variance (this is a consequence of the assumption~(H4)) and
mean 
\[
\E_z(\tilde{Z}(1)-z) ~=~ \sum_{z'\in\Z^d} (z'-z) \exp(\ol{a}\cdot(z'-z))\mu(z'-z) ~=~
\nabla\varphi(\ol{a}). 
\]
The last coordinate of $\nabla\varphi(\ol{a})$ is negative or zero because 
$\ol{a}\in\partial_- D$. Since
$\mu(z)=0$ for all $z=(x,y)$ with  
$y<-1$,  the twisted random walk $\tilde{Z}(t)$ starting at any point
$z\in\Z^{d-1}\times\N^*$ hits the 
hyper-plane $\Z^{d-1}\times\{0\}$ with probability $1$ and consequently, 
\be\label{e2-1}
\E_z\bigl(\exp(a\cdot Z(\tau)); \; \tau < \infty\bigr) ~=~ \E_z\bigl(\exp(\ol{a}\cdot Z(\tau)); \;
\tau < \infty\bigr)  ~=~ \exp(\ol{a}\cdot z)
\ee
for every $z\in\Z^{d-1}\times\N^*$. Finally, by  Markov
property,  for $z\in\Z^{d-1}\times\{0\}$ we get  
\begin{align*}
\E_z(\exp(a\cdot Z(\tau));
\;\tau < \infty) &=  
\sum_{z'\in\Z^{d-1}\times\N^*} p(z,z') \E_{z'}(\exp(a\cdot Z(\tau)); \; \tau < \infty)\\
&\hspace{3.5cm} + \sum_{z'\in\Z^{d-1}\times\{0\}}
p(z,z') \exp(a\cdot z') \\ 
&=~ \varphi_0(\ol{a}) \exp(\ol{a}\cdot z).
\end{align*}
The last relation is a consequence of  \eqref{e2-1}  and the equality $\ol{a}\cdot z' = a \cdot z'$ for 
$z'\in\Z^{d-1}\times\{0\}$. 
\end{proof}
By strong Markov property Lemma~\ref{lem2-2} implies that 
\begin{cor}\label{cor2-1} Under the hypotheses (H0)-(H1), for all $a\in D$, 
  $z\in\Z^{d-1}\times\N$, 
\[
\sum_{w\in\Z^{d-1}\times\{0\}} G(z,w)\exp(a\cdot w) ~=~ \begin{cases}
      (1-\varphi_0(\ol{a}))^{-1} \exp(\ol{a}\cdot z)  &\text{if $\varphi_0(\ol{a}) < 1$,}\\
+\infty &\text{if $\varphi_0(\ol{a}) \geq 1$.}
\end{cases}
\]
\end{cor}
The last statement together with Lemma~\ref{lem2-1} implies the following estimate
for Green's function.
\begin{cor}\label{cor2-2} Under the hypotheses (H0)-(H2), for any $a\in D$ such that
  $\varphi(\ol{a}) < 1$,  $z\in\Z^{d-1}\times\N$ 
 and $z'\in\Z^{d-1}\times\N^*$, 
\be\label{e2-2}
G(z,z')/G_S(0,0) ~\leq~ \exp(a\cdot(z-z'))  + \varphi_0(a)(1-\varphi_0(\ol{a}))^{-1}
\exp(\ol{a}\cdot z - a\cdot z') 
\ee
\end{cor} 
\begin{proof} Indeed, let $a\in D$ be such that $\varphi_0(\ol{a})< 1$. Then  for $z\in\Z^{d-1}\times\N$ and
$z'\in\Z^{d-1}\times\N^*$, from 
the renewal equation 
\[
G(z,z') ~=~ G_+(z,z') ~+~ \sum_{w\in\Z^{d-1}\times\{0\}} \;\, \sum_{w'\in\Z^{d-1}\times\N^*}
G(z,w) \mu_0(w,w') G_+(w',z') 
\]
combined with Lemma~\ref{lem2-1} it follows that 
\begin{align*}
\frac{G(z,z')}{G_S(0,0)} &\leq \exp(a\cdot (z-z')) ~+
\sum_{\substack{w\in\Z^{d-1}\times\{0\},\\ w'\in\Z^{d-1}\times\N^*}}
G(z,w) \mu_0(w'-w) \exp(a\cdot(w'-z'))\\
&\leq \exp(a\cdot (z-z')) ~+~ \varphi_0(a) \!\!
\sum_{w\in\Z^{d-1}\times\{0\}}
G(z,w) \exp(a\cdot(w-z')) 
\end{align*}
and hence, using by Corollary~\ref{cor2-1} we get \eqref{e2-2}
\end{proof}

\noindent
We will need moreover the following consequence of Lemma~\ref{lem2-2}. 
\begin{cor}\label{cor2-3} Under the hypotheses (H0)-(H4), every point of the set  
  $\partial_-D\setminus\partial_0 D$ has a neighborhood in which the function 
$
a \to \E_0\bigl(\exp(a\cdot Z(\tau)); \, \tau < \infty\bigr) 
$
is finite. 
\end{cor}

Now we obtain an explicit representation of the normal cone $V(a)$ for
$a\in\hat{D}$. Recall that $V(a)$ is
the normal cone to the convex set $\hat{D}=\{a\in D : \varphi_0(\ol{a})\leq 1\} = (\Theta\times\R)\cap D$ at the point $a\in\hat{D}$. If the point
$a$ belongs to the interior of the set $D$ then clearly $V(a)=\{0\}$. It is sufficient
therefore to consider the points on the boundary $
\partial\hat{D}$ of $\hat{D}$. According to the definition of the set $\hat{D}$, a point
$a$ belongs to the boundary $\partial\hat{D}$ if and only if  $\max\{\varphi(a),\varphi_0(\ol{a})\} =1$.

\begin{lemma}\label{lem2-3} Under the hypotheses (H0)-(H4), for every
  $a\in\partial\hat{D}$, 
\begin{align}\label{e2-3}
V(a) ~=~\begin{cases}\bigl\{ c_1\nabla\varphi(a) + c_2 (\nabla\varphi_0(\ol{a}) + \kappa_{\ol{a}}
\nabla\varphi(\ol{a})) : c_i
\geq 0\bigr\} &\text{if 
    $\varphi(a)=\varphi_0(\ol{a}) = 1$ }\\
&\text{ and $a\not\in\partial_0D$}\\
\bigl\{ c (\nabla\varphi_0(\ol{a}) + \kappa_{\ol{a}} \nabla\varphi(\ol{a})) :~ c\geq 0\bigr\}
&\text{if $\varphi(a) < \varphi_0(\ol{a}) = 1$} \\ 
&\text{ and $a\not\in\partial_0D$}\\
 \bigl\{ c\nabla\varphi(a) :~c \geq 0\bigr\} &\text{if either $a\in\partial_0 D$}\\
&\text{ or $\varphi_0(\ol{a})  < \varphi(a)=1$}
\end{cases}
\end{align}
with 
\be\label{e2-4}
\kappa_a = - \left.\frac{\partial\varphi_0(\a,\beta)}{\partial\beta}
\left(\frac{\partial\varphi(\a,\beta)}{\partial\beta}\right)^{-1}\right|_{(\a,\beta) = a}.
\ee
\end{lemma}
\begin{proof}
Indeed, under the hypotheses (H0)-(H4), the set $$
 D \cap D_0 = \{a\in\R^d :~ \max\{\varphi(a),\varphi_0(a)\}\leq
1\} $$ 
has a non-empty interior because $\varphi(0) = \varphi_0(0) = 1$ and 
\[
\frac{\nabla \varphi(0)}{|\nabla \varphi(0)|} + \frac{\nabla\varphi_0(0)}{|\nabla
  \varphi_0(0)|} ~=~ \frac{m}{|m|} + \frac{m_0}{|m_0|}~\not=~ 0. 
\]
Since $D \cap D_0 ~\subset~ \hat{D} = (\Theta\times\R)\cap D$, the set $\hat{D}$ has  also a non-empty
interior and by Corollary~23.8.1 of Rockafellar~\cite{R}, 
\be\label{e2-5}
V(a) ~=~ V_D(a) + V_{\Theta\times\R}(a),  \quad \forall a\in\hat{D} 
\ee
 where $V_{\Theta\times\R}(a)\subset\R^{d-1}\times\{0\}$ is the normal cone to the cylinder $\Theta\times\R$
at the point $a$ and 
\be\label{e2-6}
V_D(a) ~=~ \begin{cases} \{c\nabla\varphi(a) :~ c\geq 0\} &\text{if $\varphi(a)=1$}\\
\{0\}&\text{if $\varphi(a) <1$}
\end{cases}
\ee
is a normal cone to the set $D$ at the point $a$. Furthermore, recall that  $\Theta\times\{0\}$ is the
orthogonal projection of the set $D\cap D_0$ onto the hyper-plane
$\R^{d-1}\times\{0\}$. since the orthogonal projection onto the hyper-plane
$\R^{d-1}\times\{0\}$ of the point $a\in D$ is the same
as the orthogonal projection  of the point $\ol{a}$ from this it follows that 
\[
V_{\Theta\times\R}(a) ~=~ V_{\Theta\times\R}(\ol{a}) ~=~ V_{D\cap D_0}(\ol{a}) \cap
\left(\R^{d-1}\times\{0\}\right)  \quad \quad \forall a\in \hat{D}  
\]
where $V_{D\cap D_0}(\ol{a})$ denotes the normal cone to the set $D\cap D_0$ at the point
$\ol{a}\in\partial_-D$. Moreover, since $\varphi(\ol{a})=1$, using
again Corollary~23.8.1 of Rockafellar~\cite{R} we get  
\[
V_{D\cap D_0}(\ol{a}) ~=~ V_D(\ol{a}) + V_{D_0}(\ol{a}) ~=~ \begin{cases} \{c_1\nabla\varphi(\ol{a}) +
c_2\nabla\varphi_0(\ol{a}):~ c_i\geq 0\} &\text{if
    $\varphi_0(\ol{a})=1$,}\\
\{c\nabla\varphi(\ol{a}) :~ c\geq 0\} &\text{if
    $\varphi_0(\ol{a}) < 1$} 
\end{cases}
\] 
and hence, for any $a\in\hat{D}$, 
\[
V_{\Theta\times\R}(a) ~=~ \begin{cases} \bigl\{ c (\nabla\varphi_0(\ol{a}) +
  \kappa_{\ol{a}} \nabla\varphi(\ol{a})) :~ c\geq 0\bigr\} &\text{if
    $\varphi_0(\ol{a})=1$ and $\ol{a}\not\in\partial_0D$,}\\
 \bigl\{ c \nabla\varphi(\ol{a})  :~ c\geq 0\bigr\} &\text{if $\ol{a}\in\partial_0D$,}\\ 
\{0\} &\text{if
    $\varphi_0(\ol{a}) < 1$}.
\end{cases} 
\]
Finally, if  $a\in\hat{D}$ and $\ol{a}\in\partial_0D$ then clearly $a=\ol{a}$, and
consequently, the last relation combined with \eqref{e2-5} and \eqref{e2-6} prove
\eqref{e2-3}. 
\end{proof}

The next Lemma is needed to show that the mapping $q\to \hat{a}(q)$ is well defined. 

\begin{lemma}\label{lem2-4} Under the hypotheses (H0)-(H4), the set $$\Theta ~\dot=~ \{ \a
  \in\R^{d-1} :~ \inf_{\beta\in\R} \max\{\varphi(\a,\beta), 
\varphi_0(\a,\beta)\} \leq 1\}$$ is strictly
  convex~:~ for any two different
  points $\a,\a'\in\Theta$ and any $0<\theta<1$, the 
  point $\a_\theta = \theta\a + (1-\theta)\a'$ belongs to the interior of the set
  $\Theta$. 
\end{lemma}
\begin{proof} The set $\{a\in D
  : \varphi_0(a)\leq 1\} = \{a\in\R^d : \varphi(a) \leq 1 \; \text{and} \; \varphi_0(a) \leq
  1\}$ is compact and convex because the functions $\varphi_0$ and $\varphi$ are
  continuous and convex on $\R^d$. The
  set $\Theta$ is  therefore also compact and convex because $\Theta\times\{0\}$ is an
  orthogonal projection of the  set $\{a\in D 
  : \varphi_0(a)\leq 1\}$ on the hyper-plane $\R^{d-1}\times\{0\}$. Furthermore, remark that 
$
\Theta\subset \{\a\in\R^{d-1} ~:~ \inf_\beta \varphi(\a,\beta) \leq 1\} 
$ 
and that  the mapping $a=(\a,\beta)\to \a$ determines a
  homeomorphism between the set $\partial_- D$ and the set $\{\a \in\R^{d-1}, \; \inf_\beta
  \varphi(\a,\beta) \leq 1\}$.  Let $\a\to (\a,\beta_\a)$ denote the inverse mapping to such a homeomorphism. 
Since for every $\alpha\in\R^{d-1}$, the
  function $\beta \to \varphi_0(\a,\beta)$ is increasing then a point  $\a\in\R^{d-1}$
  satisfying the inequality $\inf_\beta
  \varphi(\a,\beta) \leq 1$ belongs to the set $\Theta$ if and only if
  $\varphi_0(\a,\beta_\a) \leq 1$ and consequently, 
\[
\Theta ~=~  \{\a\in\R^{d-1} ~:~ \inf_\beta \varphi(\a,\beta) \leq 1 \; \text{and} \; \varphi_0(\a,\beta_\a) \leq 1\}. 
\]
Under the hypotheses (H2),  the set $D$ is strictly
convex because the function $\varphi$ is strictly convex. The set $
\{\a\in\R^{d-1} :~ \inf_\beta
\varphi(\a,\beta) \leq 1\}$  is therefore also strictly convex and hence, to prove that
the set $\Theta$ is strictly 
  convex it is sufficient to show that the function
$\a \to \varphi_0(\a,\beta_\a)$ is strictly convex on  $\{\a\in\R^{d-1} : \inf_\beta
\varphi(\a,\beta) < 1\}$. For this we use Lemma~\ref{lem2-2}. Recall that by
Lemma~\ref{lem2-2}, for $a=(\a,\beta_\a)\in\partial_- D$,  
\[
\varphi_0(\a,\beta_\a) ~=~ \E_0\bigl(\exp(a\cdot Z(\tau)); \, \tau < \infty\bigr) ~=~
\sum_{x\in\Z^{d-1}} \P_0(X(\tau) =x) ~\exp(\a\cdot x) 
\]
where $\tau$ is the first time when the process $Z(t)=(X(t),Y(t))$ returns to the boundary
hyper-plane $\Z^{d-1}\times\{0\}$. 
Since under the hypotheses of our lemma, the
function $\varphi_0$ is finite  everywhere on $\R^d$ then the series at the
right hand side of the above relation converge on $\{\a\in\R^{d-1} ~:~ \inf_\beta
\varphi(\a,\beta) \leq 1\}$. By dominated convergence theorem, from this it follows that
the function $\a\to\varphi_0(\a,\beta_\a)$ is infinitely differentiable on  $\{\a\in\R^{d-1} ~:~ \inf_\beta
\varphi(\a,\beta) < 1\}$ and that its  Hessian matrix 
\[
Q(\a) ~=~ \left(\frac{\partial^2  \varphi_0(\a,\beta_\a)
 }{\partial\alpha_i\partial\alpha_j}\right)_{1\leq i,j\leq d-1} 
\]
satisfies the equality  
\[
\xi\cdot Q(\a)\xi ~=~ \sum_{x\in\Z^{d-1}} ~e^{\a\cdot x} ~(\xi\cdot x)^2 \P_0(X(\tau) =x) 
\]
for any $\xi\in\R^{d-1}$ whenever 
$\inf_\beta \varphi(\a,\beta) < 1$. Since the Markov
process $Z(t)$ is irreducible, then for every non-zero vector $\xi\in\R^{d-1}$ there is
$x\in\Z^{d-1}$ such that $ (\xi\cdot x)^2 \P_0(X(\tau) =x)  ~>~ 0$
and consequently, $\xi\cdot Q(\a)\xi > 0$. 
 This proves that the function $\a\to \varphi(\a,\beta_\a)$ is strictly convex on  the set 
$\{\a\in\R^{d-1} :~ \inf_\beta \varphi(\a,\beta) < 1\}$.  Lemma~\ref{lem2-4} is 
 proved. 
\end{proof}

We are ready now to get the following statement.

\begin{lemma}\label{lem2-5} Under the hypotheses (H0)-(H4), for every non-zero vector
  $q\in\R^{d-1}\times\R_+$, there 
  is a unique point $\hat{a}(q)\in\hat{D}\cap \partial_+ D$ for which $q\in V(\hat{a}(q))$. 
\end{lemma}
\begin{proof} Recall that for $\hat{a}\in\hat{D}$, the vector $q$
  belongs to the normal cone $V(\hat{a})$ to the set $\hat{D}$ if and only if 
\be\label{e2-7}
\sup_{a\in \hat{D}} ~a\cdot q ~=~ \hat{a}\cdot q.
\ee
Since under the hypotheses (H0)-(H4), the set $\hat{D}$ is compact and
  non-empty, for every $q\in{\cal S}^d$ there is $\hat{a}=\hat{a}(q)\in
  \hat{D}$ for which this equality holds. It is clear that for 
 $q\not=0$,  such a point $\hat{a}(q)$ belongs to the boundary
  $\partial\hat{D}$ of the set 
$\hat{D}$. Moreover, Lemma~\ref{lem2-3} shows that  for $q\in\R^{d-1}\times]0,+\infty[$, 
\[ 
\hat{a}(q) \in \hat{D}\cap \partial_+ D.
\]
For $q\in\R^{d-1}\times\{0\}$,  a point $\hat{a} = \hat{a}(q)$ satisfying the equality \eqref{e2-7}
can be non-unique~: if the equality \eqref{e2-7} holds for some
$\hat{a}\in\partial(\hat{D})$  then 
\[
\sup_{a\in \hat{D}} ~a\cdot q ~=~ \tilde{a}\cdot q
\]
 for all $\tilde{a}\in \partial\hat{D}$ having the same first $d-1$ coordinates as the point $\hat{a}$.
Remark however that for every $a\in \partial\hat{D}$,  there is a unique point $\hat{a}\in\hat{D}\cap
  \partial_+ D$
with the same first $d-1$ coordinates as the point $a$ and hence without any 
restriction of generality we can assume that $\hat{a}(q)\in \hat{D}\cap
\partial_+D$. 

We have shown that for every non-zero vector
  $q\in\R^{d-1}\times\R_+$, there
  is a point $\hat{a}(q)\in\hat{D}\cap \partial_+ D$ for which $q\in 
  V(\hat{a}(q))$. To complete the proof of our lemma it is now sufficient to
  show that such a point  is unique. Suppose that there are two different points
  $\hat{a}(q),\tilde{a}(q)\in\hat{D}\cap \partial_+ D$ for which $q\in
  V(\hat{a}(q))\cap V(\tilde{a}(q))$. Then 
\[
\sup_{a\in \hat{D}} ~a\cdot q ~=~ \tilde{a}(q)\cdot q ~=~ \hat{a}(q)\cdot q
\] 
 and consequently,  for every
$\theta\in[0,1]$,  
\[
\sup_{a\in \hat{D}} ~a\cdot q ~=~ \bigl(\theta \tilde{a}(q) + (1-\theta) \hat{a}(q)\bigr)\cdot q.
\] 
The last equality shows  that the point $a_\theta \dot=~ \theta \tilde{a}(q) +
(1-\theta) \hat{a}(q)$ belongs to 
the boundary of the set $\hat{D}$ and that $
q\in V(a_\theta)$. 
Recall now that under the hypotheses of our lemma, the set $D$ is strictly
convex and consequently, for $0<\theta<1$, the point $a_\theta ~\dot=~ \theta \tilde{a}(q) + (1-\theta)
\hat{a}(q)$ belongs to the interior of the set $D$. Hence, the normal
cone $V_D(a_\theta)$ to the set $D$ at the point $a_\theta$ is zero and the normal cone
$V(a_\theta)$ to $\hat{D}$ at the point $a_\theta$ coincide with the 
normal cone $V_{\Theta\times\R}(a_\theta)$ to the set $\Theta\times\R$ at $a_\theta$ (this
is a consequence of Corollary~23.8.1 of \cite{R}). From this it follows that  
\be\label{e2-8}
q\in V(a_\theta) ~=~ V_{\Theta\times\R}(a_\theta) \subset\R^{d-1}\times\{0\}.
\ee
For $q\in\R^{d-1}\times]0,+\infty[$, the point $
\hat{a}(q)=\tilde{a}(q)$ is therefore unique. For $q\in\R^{d-1}\times\{0\}$, 
  \eqref{e2-8} shows that the first $d-1$ coordinates of the points $a_\theta$ and
  $a_{\theta'}$ are the same for all $0<\theta<\theta' < 1$ because by Lemma~\ref{lem2-4},  the set
 $\Theta$ is strictly convex. Letting $\theta\to 0$
and $\theta'\to 1$ we conclude  that the first $d-1$ coordinates of the points 
$\hat{a}(q))$ and $\tilde{a}(q)$ are the same.  This proves that  
$\hat{a}(q)=\tilde{a}(q)$ because $\hat{a}(q), \tilde{a}(q) \in\partial_+ D$ and 
the orthogonal projection determines a one to one mapping from $\partial_+ D$ to $\R^{d-1}\times\{0\}$. 
\end{proof}

Lemma~\ref{lem2-3} and Lemma~\ref{lem2-5} imply the following statement.

\begin{cor}\label{cor2-4} Under the hypotheses (H0)-(H4), for every $q\in
  \R^{d-1}\times]0,+\infty[$, the following assertions hold :~
\begin{enumerate}
\item[1)] there is a unique vector $\gamma_q\in \R^{d-1}\times\{0\}$
      for which the vector $q-\gamma_q$ belongs to the normal cone to the set $D$ at the
      point $\hat{a}(q)$ and  $\gamma_q, \, q-\gamma_q \in V(\hat{a}(q))$.
\item[2)] $\varphi_0\left(\ol{\hat{a}(q)}\right) = 1$ whenever $\gamma_q \not= 0$.  
\end{enumerate}
\end{cor}

\section{Harmonic functions}\label{sec-hf} The
harmonic function of the Markov process $(Z(t))$ are now identified. The main result of this section is the
following proposition.

\begin{prop}\label{pr3-1}  Under the hypotheses (H0)-(H4), the following assertions hold. 

\noindent 
1) A non-negative  function $h$ is harmonic for the Markov process $(Z(t))$ if and only if
 there  is a positive measure $\nu_h$ on $\hat{D}\cap\partial_+D = (\Theta\times\R)\cap\partial_+ D$ such that 
\be\label{e3-1}
h(z) ~=~ \int_{(\Theta\times\R)\cap\partial_+ D} h_{a}(z) \, d\nu_h(a), \quad \quad \forall
z\in\N^*\times\Z^{d-1}.  
\ee
2) For every $a=(\a,\beta)\in(\Theta\times\R)\cap\partial_+ D$ with $\a\in\R^{d-1}$ and $\beta\in\R$, the
constant multiples of the function $h_{a}$ defined by \eqref{e1-8} are the 
only non-negative harmonic functions for which 
\be\label{e3-2}
\sup_{x\in\R^{d-1}} \exp(- \a\cdot x)  h(x,y) ~<~ +\infty, \quad \quad \forall
y\in\N. 
\ee
3) The  constant multiples of the functions $h_{a}$ with
  $a\in(\Theta\times\R)\cap\partial_+ D$, are the only minimal harmonic  
functions of the Markov process $(Z(t))$. 
\end{prop}

In order to prove this result we use the properties of Markov-additive processes. Recall
that a Markov process 
$(A(t),M(t))$ on a countable set $\Z^{d}\times E$ with transition probabilities  
  $p\bigl((x,y),(x',y')\bigr)$ is called {\em Markov-additive} if 
\[
p\bigl((x,y),(x',y')\bigr) ~=~ p\bigl((0,y),(x'-x,y')\bigr)
\]
for all $x,x'\in\Z^{d}$, $y,y'\in E$. The first component $A(t)$ is an {\em additive} part of the process
$(A(t),M(t))$,  and $M(t)$ is its {\em Markovian part}. 

According to this definition, the Markov process $Z(t)=(X(t),Y(t))$ is Markov-additive
with an additive part $X(t)$ taking the values in $\Z^{d-1}$ and Markovian part $Y(t)$
taking the values in $\N$. Under the hypotheses (H1), its {\it Feynman-Kac transform}
matrix ${\cal P}(\alpha) = \bigl( {\cal P}(\alpha, y,y'), \; y,y'\in \N\bigr)$
with $\a\in\R^{d-1}$ and 
\[
{\cal P}(\alpha, y,y') ~=~ \E_{(0,y)}\bigl( \exp(\alpha\cdot X(1)) ; \, Y(1) = y'\bigr)
\] 
is irreducible and  the limit
\[
\lambda(\alpha) ~=~ \limsup_n \frac{1}{n} \log {\cal P}^{(n)}(\alpha, y,y') 
\]
does not depend on $y,y'\in \N$ (see~\cite{Seneta}). The quantity $e^{\lambda(\alpha)}$ is usually called
  {\it spectral radius} and $e^{-\lambda(\alpha)}$ is the {\it convergence parameter} of the
  transform matrix ${\cal P}(\alpha)$. By Proposition~3.1 of
  Ignatiouk~\cite{Ignatiouk:06}, every non-zero minimal harmonic 
  function $h$ of the Markov process $Z(t)$ is of the form 
\be\label{e3-3}
h(x,y) = \exp(\a\cdot x) h(0,y), \quad \quad \forall (x,y)\in\Z^{d-1}\times\N. 
\ee 
with some  $\alpha\in\R^{d-1}$  satisfying
  the inequality $\lambda(\alpha) \leq 0$.   The following lemma identifies 
   the function $\alpha\to \lambda(\alpha)$. 

\begin{lemma}\label{lem3-1} Under the hypotheses $(H0)-(H4)$, 
\be\label{e3-4}
\lambda(\alpha) ~=~ \inf_{\beta\in\R} \log \max\{\varphi(\alpha,\beta),
\varphi_0(\alpha,\beta)\} , \quad \quad \forall
\a\in\R^{d-1}.
\ee
\end{lemma}
\begin{proof} Remark first of all that for any $(\a,\beta)\in\R^{d-1}\times\R$, the
  exponential function $f(y) = \exp(\beta y)$ on $\N$ satisfies the inequality 
\[
{\cal P}(\alpha) f (y) ~=~ \E_{(0,y)}(\exp(\a\cdot X(1) + \beta Y(1))) \leq \max\{\varphi(\alpha,\beta),
\varphi_0(\alpha,\beta)\} f(y) \quad \forall y\in\N. 
\]
From this it follows that  $
\lambda(\a) \leq \log \max\{\varphi(\alpha,\beta),
\varphi_0(\alpha,\beta)\}$ for all $(\a,\beta)\in\R^{d-1}\times\R$ 
(see Seneta~\cite{Seneta} for more details)  and consequently, 
\[
\lambda(\alpha) ~\leq~ \inf_{\beta\in\R} \log \max\{\varphi(\alpha,\beta),
\varphi_0(\alpha,\beta)\} , \quad \quad \forall
\a\in\R^{d-1}.
\]
Furthermore, let $\tau$ denote the first time when the process $(Z(t))$ hits the
boundary hyperplane $Z^{d-1}\times\{0\}$. Then for $y,y' > 0$, $y,y'\in\N$, 
\begin{align*}
\lambda(\a) &~=~ \limsup_n \frac{1}{n} \log \E_{(0,y)}\bigl( \exp(\alpha\cdot X(n)) ; \,
Y(n) = y'\bigr) \\
&~\geq~ \limsup_n \frac{1}{n} \log \E_{(0,y)}\bigl( \exp(\alpha\cdot X(n)) ; \, Y(n) =
y', \; \tau > n\bigr) ~=~ \inf_{\beta\in\R} \log \varphi(\a,\beta)\\ 
\end{align*}
where the last relation is proved by 
Lemma~5.1 of Ignatiouk~\cite{Ignatiouk:06}. For those $\a\in \R^{d-1}$ for which the right
hand side of \eqref{e3-4}  is equal to right hand side of the last relation, the equality
\eqref{e3-4} is therefore verified. Suppose now that 
\[
\inf_{\beta\in\R}  \max\{\varphi(\alpha,\beta),
\varphi_0(\alpha,\beta)\} ~>~ \inf_{\beta\in\R} \varphi(\a,\beta). 
\]
In this case, the minimum of the function $\beta \to \max\{\varphi(\alpha,\beta),
\varphi_0(\alpha,\beta)\}$ is achieved at a point $\hat{\beta}_\a\in\R$  where  
\[
\varphi(\alpha,\hat{\beta}_\a) = \varphi_0(\alpha,\hat{\beta}_\a) \quad \quad \text{and}
\quad \quad \frac{\partial}{\partial\beta}\varphi(\a, \hat{\beta}_\a) < 0. 
\]
Under the hypotheses $(H0)-(H4)$, the twisted Markov chain $(\tilde{Y}(t))$ on $\N$
having transition probabilities  $
\tilde{p}(y,y') ~=~ {\cal P}(\a, y, y') \exp(\hat{\beta}_\a(y'-y))/\varphi(\alpha,\hat{\beta}_\a)$ 
is irreducible and satisfies the conditions of Foster's criterion of positive recurrence
(see Corollary~8.7 in~\cite{PhilippeRobert}) with the test function
$f(y)=y$ : 
\[
\E_0(\tilde{Y}(1)) ~<~ +\infty \quad \text{and} \quad \E_y(\tilde{Y}(1))   ~=~ y +
\frac{\partial}{\partial\beta}\varphi(\a, \hat{\beta}_\a) ~<~ y, \quad \forall y>0.
\]
The  Markov chain  $(\tilde{Y}(t))$ is therefore positive recurrent and consequently,
\[
\limsup_{n\to\infty} \frac{1}{n} \log \tilde{p}^{(n)}(y,y') ~=~ 0, \quad \quad  \forall
y,y'\in\N. 
\]
The last relation together with the equality $
{\cal P}^{(n)}(\a, y, y) =
\tilde{p}^{(n)}(y,y) (\varphi(\a,\hat{\beta}_\a))^n$
shows that $\lambda(\a) = \log
\varphi(\a,\hat{\beta}_\a)$ from which it follows \eqref{e3-4}.   
\end{proof}

Lemma~\ref{lem3-1} proves that $\lambda(\alpha) \leq 0$ if and only if
$\alpha\in\Theta$ and hence, using Proposition~3.1 of
  Ignatiouk~\cite{Ignatiouk:06} we get   
\begin{cor}\label{cor3-1} Under the hypotheses $(H0)-(H4)$,   every minimal harmonic function
  $h$ of the Markov process $(Z(t))$ satisfies the equality   \eqref{e3-3} 
with some $\a\in \Theta$. 
\end{cor}

Now we identify the minimal harmonic functions satisfying the equality  \eqref{e3-3}.

\begin{lemma}\label{lem3-2} Under the hypotheses (H0)-(H4),  for every point
  $a=(\a,\beta)\in(\Theta\times\R)\cap  
  \partial_+ D$, the constant multiples of $h_a$ are the only minimal non-negative harmonic functions of the
  Markov process $(Z(t))$ for which the equality \eqref{e3-3} holds with a given $\alpha\in\Theta$.  
\end{lemma}
\begin{proof}   Let $a=(\a,\beta)\in(\Theta\times\R)\cap  
  \partial_+ D$. 
Straightforward calculation shows
that the function $h_a$ is non-negative and harmonic for the Markov process
$(Z(t))$. Recall that a non-zero harmonic function $h\geq 0$ is called 
  {\em minimal} if for any non-zero harmonic function $h'\geq 0$, the inequality $h'\leq h$
  implies that $h'=c h$ with some constant $c > 0$.  
To prove our Lemma it is therefore
  sufficient to show that if $h\not=0$ is a minimal non-negative harmonic functions of the
  Markov process $(Z(t))$ for which \eqref{e3-3} holds with a given
  $\alpha$ then 
\be\label{e3-5}
h \geq c h_a
\ee 
with some $c > 0$. For this we first show
  that every such a function $h\not=0$ satisfies the inequality 
\be\label{e3-6}
h(z) ~\geq~ h(0) \exp(\ol{a}\cdot z) ~>~ 0 \quad \quad \text{for all $z\in\Z^{d-1}\times\N$}.
\ee
Indeed, let $h$ be a non-zero minimal non-negative
  harmonic functions for which the equality \eqref{e3-3} holds with a given
  $\alpha$. Then $h(z) > 0$ for all $z\in\Z^{d-1}\times\N$ because the Markov
  process $Z(t)$ is irreducible. 
Moreover,  according to the definition of the mapping
  $a\to\ol{a}$,  from \eqref{e3-3}  it follows
  that 
\be\label{e3-7}
h(z) ~=~ h(0) \exp(\a\cdot x) ~=~ h(0) \exp(\ol{a}\cdot z) ~>~ 0 \quad \text{for any $z=(x,0)\in\Z^{d-1}\times\{0\}$.}
\ee
Hence, for $z\in\Z^{d-1}\times\{0\}$ the inequality \eqref{e3-6} holds with the
equality. 
Furthermore, for $\tau = \inf\{t > 0 : Y(t) = 0\}$, the sequence $h(Z(n\wedge\tau))$ is a
  martingale relative to the natural filtration and $h(Z(n\wedge\tau)) = h(0) \exp(\alpha\cdot
  X(\tau)) $ whenever $\tau\leq n$. Hence, for any
  $z=(x,y)\in\Z^{d-1}\times\N$ with $y>0$ we have 
\[
h(z) ~=~ \E_{z}(h(Z(n\wedge\tau)) ~\geq~ h(0) E_{z}\Bigl(\exp(\a\cdot
X(\tau)); \; \tau \leq n\Bigr), \quad \forall n\in\N 
\]
and consequently, letting   $n\to\infty$ and using Fatou
lemma  we obtain  
\[
h(x,y)  ~\geq~ h(0) E_{z}\Bigl(\exp(\a\cdot
X(\tau)); \; \tau < \infty\Bigr)
\]
The last inequality combined with Lemma~\ref{lem2-2} proves \eqref{e3-6} for
$z=(x,y)\in\Z^{d-1}\times\N$  with $y>0$. The inequality \eqref{e3-6} is therefore
verified. 

Recall now that $
(\Theta\times\R)\cap\partial_+D ~=~ \{a\in \partial_+D :~ \varphi_0(\ol{a})\leq 1\}$ 
where $\ol{a}$ is a point on the boundary
$\partial_-D$ having the same $d-1$ first coordinates as the point $a$. From now on the
proof of \eqref{e3-5} is different
  in each of the following cases :
\begin{itemize}
\item[--]{\em case 1} : when $\varphi_0(\ol{a}) = 1$,
\item[--]{\em case 2} : when $\varphi_0(\ol{a}) < 1$.
\end{itemize}
If $\varphi_0(\ol{a}) = 1$ then from \eqref{e1-8} it follows that $h_a(z) = \exp(\ol{a}\cdot z)$ for all $z\in\Z^{d-1}\times\N$ and hence,
the inequality \eqref{e3-6} proves \eqref{e3-5} with $c=h(0)$. For all those
$a=(\a,\beta)\in(\Theta\times\R)\cap \partial_+ D$  for which $\varphi_0(\ol{a}) = 1$,
Lemma~\ref {lem3-2} is therefore proved.  

Suppose now that $\varphi_0(\ol{a}) < 1$  and let  $ 
h_+(z) = h(z) -  h(0)\exp(\ol{a}\cdot z)$.  Then the inequality  \eqref{e3-6} shows that
the function $h_+$ is non-negative, the equality \eqref{e3-7} implies that 
\be\label{e3-8}
h_+(z) ~=~ 0, \quad \quad \text{for any $z=(x,0)\in\Z^{d-1}\times\{0\}$,}
\ee
and from the equality \eqref{e3-3} it follows that 
\be\label{e3-9}
h_+(x,y) = \exp(\a\cdot x) h_+(0,y), \quad \quad \text{for all $z=(x,y)\in\Z
^{d-1}\times\N$.} 
\ee
Moreover,  straightforward
calculations show that  for $z=(x,0)\in\Z^{d-1}\times\{0\}$, 
\be\label{e3-10}
\E_{z}(h_+(Z(1))) ~=~ (1 - \varphi_0(\ol{a})) \exp(a\cdot z)h(0) ~>~ 0, \quad \forall x\in\Z^{d-1},  
\ee
and for $z=(x,y)\in\Z^{d-1}\times\N$ with $y >0$, 
\be\label{e3-11}
\E_{z}(h_+(Z(1))) ~=~  h(z) -
\varphi(\ol{a}) \exp(\ol{a}\cdot z)h(0) ~=~ h(z) - \exp(\ol{a}\cdot z)h(0)  ~=~ h_+(z).
\ee 
According to the definition of the
Markov process $(Z(t))$, relations \eqref{e3-8} and \eqref{e3-11} show that the function
$h_+$ satisfies the equality 
\be\label{e3-12}
\sum_{z'\in\Z^{d-1}\times\N^*} \mu(z'-z) h_+(z') ~=~ h_+(z) \quad \quad \forall
z\in\Z^{d-1}\times\N^*,
\ee
and from \eqref{e3-10} it follows that $h_+ \not\equiv 0$. Under the hypotheses of our
lemma, Proposition~2.1 and Proposition~5.1 of  
Ignatiouk~\cite{Ignatiouk:06} prove that the only non-negative
non-zero functions 
satisfying the equalities  \eqref{e3-9} and \eqref{e3-12} are the constant multiples of
\[
h_{a,+}(z) ~=~ \begin{cases} \exp(a\cdot z) - \exp(\ol{a}\cdot z) &\text{if
    $a\not\in\partial_0 D$,}\\
y ~\exp(a\cdot z) &\text{if
    $a\in\partial_0 D$,} \quad \quad \;  z=(x,y)\in\Z^d\times\N^*. 
\end{cases}
\]
Hence,  $h_+(x,y) = c h_{a,+}(x,y)$ 
for all $(x,y)\in\Z^{d-1}\times \N^*$ with some $c>0$, and consequently, 
\[
h(z) = h(0)\exp(\ol{a}\cdot z) + c h_{a,+}(z) \quad \quad \quad \forall \;  z\in\Z^{d-1}\times\N^*.
\] 
To complete the proof of our lemma it is sufficient now to notice that  $h_a(z) = C
\exp(\ol{a}\cdot z)  + h_{a;+}(z)$ with 
\[
C ~=~ \begin{cases} \frac{\partial}{\partial \beta}\varphi_0(a)/(1-\varphi_0(\ol{a})
  &\text{ if $a\in\partial_0 D$}\\
(\varphi_0(a) - \varphi_0(\ol{a}))/(1-\varphi_0(\ol{a}) &\text{ otherwise} 
\end{cases} 
\]
from which it follows that  $h(z) \geq \min\{c, h(0)/C\} h_a(z)$.
\end{proof}

Lemma~\ref{lem3-2} combined with Corollary~\ref{cor3-1} implies the following statement.

\begin{cor}\label{cor3-2} Under the hypotheses (H0)-(H4), every minimal harmonic function
  of the Markov process $(Z(t))$ is of the form   $h = ch_a$ 
with some $c>0$ and $a\in(\Theta\times\R)\cap\partial_+ D$.  
\end{cor}
\begin{proof} To get this statement from Corollary~\ref{cor3-1}  and Lemma~\ref{lem3-2} it
  is sufficient to notice that the orthogonal projection onto the hyper-plane
  $\R^{d-1}\times\{0\}$ determines a
  homeomorphism from $(\Theta\times\R)\cap\partial_+ D$ to $\Theta\times\{0\}$. 
\end{proof}

\noindent{\bf Proof of Proposition~\ref{pr3-1}.} The proof of this proposition uses
Corollary~\ref{cor3-2} and the same arguments as in the proof of Proposition~5.1 of
Ignatiouk~\cite{Ignatiouk:06}. The main steps of this proof are the following.

By the Poisson-Martin representation theorem
(see Woess~\cite{Woess}), every 
non-negative harmonic function of the Markov process $(Z(t))$ is of the form 
\[
h(z) ~=~ \int_{\partial_m(\Z^{d-1}\times\N)} K(z,\gamma) \, d\tilde\nu_h(\gamma), \quad \quad
\forall z\in\Z^{d-1}\times\N^*
\]
with some Borel measure $\tilde\nu_h\geq 0$ on the minimal Martin boundary
$\partial_m(\Z^{d-1}\times\N)$.   Recall that 
 $K(z,\gamma)$ is
the Martin kernel of the Markov process $(Z(t))$, the mapping $\gamma\to
K(z,\gamma)$ is continuous on $\partial_m(\Z^{d-1}\times\N)$ for every
$z\in\Z^{d-1}\times\N$ and for
every $\gamma\in\partial_m(\Z^{d-1}\times\N)$, according to the definition of the minimal
Martin boundary (see Woess~\cite{Woess}), the function $z \to K(z,\gamma)$ 
is a minimal harmonic function for the Markov process $(Z(t))$ with $K(z_0,\gamma) =
1$. By
Corollary~\ref{cor3-2}, we have therefore 
\[
K(z,\gamma) = c_\gamma h_{a(\gamma)}(z) \quad \text{for all $z\in\Z^{d-1}\times\N$}
\]
with some $a(\gamma)=(\a(\gamma),\beta(\gamma))\in(\Theta\times\R)\cap\partial_+ D$ and $c_\gamma=
1/h_{a(\gamma)}(z_0)$. For 
$z_0=(x_0,y_0)$ and $z=(x,y_0)$, the mapping  
\[
\gamma ~\to~ K(z,\gamma) ~=~ \exp\bigl(\a(\gamma)\cdot(x-x_0)\bigr) ~>~ 0 
\]
 is therefore continuous on $\partial_m(\Z^{d-1}\times\N^*)$ for any
    $x\in\Z^{d-1}$.  This proves that the mapping $\gamma\to\a(\gamma)$ from
$\partial_m(\Z^{d-1}\times\N^*)$ to $\Theta$ is  continuous. 
The mapping
$\gamma\to a(\gamma)\in(\Theta\times\R)\cap\partial_+ D$ is therefore also continuous on
  $\partial_m(\Z^{d-1}\times\N^*)$ because the mapping $(\alpha,\beta) \to \alpha$ defines
a homeomorphism from $(\Theta\times\R)\cap\partial_+ D$  to  $\Theta$.  From
this it follows that  the integral representation \eqref{e3-1} holds with the positive
Borel measure $\nu_h$ on 
$(\Theta\times\R)\cap \partial_+ D$ defined by
\[
\nu_h(B) ~=~ \int_{\{\gamma : a(\gamma)\in B\}} c_\gamma \, d\tilde\nu_h(\gamma) 
\]
for every Borel subset $B\subset(\Theta\times\R)\cap\partial_+ D$.  The first assertion of
Proposition~\ref{pr3-1} is therefore proved.

To prove the second assertion it is sufficient to show that a non-zero harmonic
function $h\geq 0$ satisfies  \eqref{e3-2} with some $a=(\a,\beta)\in(\Theta\times\R)\cap\partial_+ D$ if and
only if $$\supp(\nu_h) =\{a\}.$$ For every 
$a\in(\Theta\times\R)\cap\partial_+ D$, the function $h_{a}$  satisfies  \eqref{e3-2} and 
is harmonic for the Markov process $(Z(t))$. 
Conversely, if  $\supp(\nu_h) \not=\{\hat{a}\}$ for some
$\hat{a}=(\hat\a,\hat\beta)\in(\Theta\times\R)\cap\partial_+ D$, 
then there is an open ball $B(a_0,\eps)$
in $\R^d$  centered at some point  $a_0=(\a_0,\beta_0)\in
(\Theta\times\R)\cap\partial_+ D$ and having a radius  $\eps >0$ such that 
$\nu_h(B(a_0,\eps)\cap(\Theta\times\R)\cap\partial_+ D) >0$ and there is 
$x_0\in\Z^{d-1}$ such that $(\a
- \hat\a)\cdot x_0 >0$  for all $a=(\a,\beta)\!\in\!B(a_0,\eps)$.
Using the integral
representation \eqref{e3-1} and 
 Fatou lemma from this it follows that  
\begin{align*}
\sup_{x\in\R^{d-1}} e^{-\hat{\a}\cdot x} h(x,y)  &\geq~ \limsup_{n\to\infty} ~e^{- n\hat{\a}\cdot
  x_0} h(n x_0,y) \\ &\geq~   \limsup_{n\to\infty} ~\int_{B(a_0,\eps)\cap(\Theta\times\R)\cap\partial_+
  D} e^{n (\a - \hat\a)\cdot x_0} h_{a,+}(0,y)  \, d\nu_h(a) \\ &\geq~
  ~\int_{B(a_0,\eps)\cap(\Theta\times\R)\cap\partial_+ 
  D} \lim_{n\to\infty} e^{n (\a - \hat\a)\cdot x_0} h_{a,+}(0,y)  \, d\nu_h(a)   =
 + \infty.
\end{align*}
The second assertion of Proposition~\ref{pr3-1} is proved. 

Finally,   if a non-negative harmonic function $h$ satisfies 
the inequality $h \leq h_{a}$ for some $a\in(\Theta\times\R)\cap\partial_+ D$ then for $h$ the inequality
\eqref{e3-2} holds with the same $a$ and
consequently $h=ch_{a}$ for some $c\geq 0$. For every
$a\in(\Theta\times\R)\cap\partial_+D$, the harmonic
function $h_{a} > 0$ is therefore minimal and  conversely, by Corollary~\ref{cor3-2}, every
minimal harmonic function of the Markov process $(Z(t))$ is of the form $c h_{a}$ with
some $c>0$ and $a\in(\Theta\times\R)\cap\partial_+ D$. Proposition~\ref{pr3-1} is proved.

\section{Communication condition}\label{sec-cc}

\noindent
{\bf Definition : } {\em A discrete time Markov chain $(Z(t))$ on $\Z^d$ is said to satisfy
  communication condition on $E\subset\Z^d$ if there exist 
$\theta >0$ and $C>0$ such that for any $z,z'\in E$ there is a sequence of
  points $z_0, z_1, \ldots,z_n\in E$  with $z_0=z$, $z_n=z'$ and 
    $n\leq C|z'-z|$ such that  $|z_i-z_{i-1}| \leq C$ and $\P_{z_{i-1}}(Z(1) = z_i) \geq
 \theta$ for all $ i=1,\ldots,n$. 
}

\begin{prop}\label{pr4-1} Under the hypotheses (H0)-(H3), the  Markov process  $(Z(t))$
  satisfies communication condition on $\Z^{d-1}\times\N$. 
\end{prop}
\begin{proof} Recall that on the half-space $\Z^{d-1}\times\N^*$,
  the Markov process $(Z(t))$ behaves as a 
homogeneous random walk $(S(t))$ on $\Z^d$ having transition probabilities
$p(z,z')=\mu(z'-z)$. Let $(Z_+(t))$ denote a sub-stochastic random walk on
$\Z^{d-1}\times\N^*$ with transition matrix $
\left(p(z,z')=\mu(z'-z), \; z,z'\in\Z^{d-1}\times\N^*\right)$.  
Such a Markov process is identical to the homogeneous random walk $(S(t))$ until the first time
when $(S(t))$ hits the boundary hyperplane $\Z^{d-1}\times\{0\}$ and dies when $(S(t))$ hits
$\Z^{d-1}\times\{0\}$.  By Lemma~4.1 of Ignatiouk~\cite{Ignatiouk:06}, the Markov process
$(Z_+(t))$ satisfies communication condition on $\Z^{d-1}\times\N^*$ ~: there exist  
$\theta >0$ and $C>0$ such that for any $z,z'\in\Z^{d-1}\times\N^*$ there is a sequence of
  points $z_0, z_1, \ldots,z_n\in\Z^{d-1}\times\N^*$  with $z_0=z$, $z_n=z'$ and 
    $n\leq C|z'-z|$ such that  
\[
|z_i-z_{i-1}| \leq C \quad \text{ and } \quad 
\mu(z_i-z_{i-1}) \geq \theta,  \quad \forall \;
 i=1,\ldots,n. 
\]   
Since the Markov
process $(Z(t))$ has the same transition probabilities on the set $\Z^{d-1}\times\N^*$ as
$(Z_+(t))$, we conclude that $(Z(t))$ also satisfies communication condition on $\Z^{d-1}\times\N^*$ with the
same constants $C>0$ and $\theta > 0$. 
Moreover, the Markov process $(Z(t))$ is
irreducible and its transition probabilities are invariant with respect to the shifts on
$z\in\Z^{d-1}\times\{0\}$. Hence, 
there are $w,w'\in\Z^{d-1}\times\N^*$ 
such that $
p(z,z+w) = \mu_0(w) > 0$ and $p(z + w',z) = \mu(-w') >0$
for all $z\in\Z^{d-1}\times\{0\}$. From this it follows that the Markov process $(Z(t))$
satisfies  communication condition  on $\Z^{d-1}\times\N$ with another
constants $C'= C + |w| + |w'|$ and $\theta' =
\min\{\theta, \mu_0(w), \mu(-w')\}$ ~:~ for any
$z,z'\in\Z^{d-1}\times\N$ there is a sequence of 
  points $z_0, z_1, \ldots,z_n\in\Z^{d-1}\times\N$  with $z_0=z$, $z_n=z'$  and 
    $n\leq C'|z'-z|$ such that  $
 |z_i-z_{i-1}| \leq C'$ and $\P_{z_{i-1}}(Z(1) = z_i) ~\geq~ \theta'$ for all
 $i=1,\ldots,n$ where $z_1 = z+w$ if $z\in\Z^{d-1}\times\{0\}$ and $z_{n-1} = z' -w'$ if
  $z'\in\Z^{d-1}\times\{0\}$. 
\end{proof}

\section{Large deviation estimates}\label{sec-ld} 
In this section, we obtain large deviation estimates for 
Green's function of the Markov processes $(Z(t))$ and $(Z_+(t))$ by using 
sample path large deviation  properties of scaled processes 
$Z^\eps(t)~=~ \eps Z([t/\eps])$ and $Z_+^\eps(t)= \eps Z_+([t/\eps])$. Recall that $(Z_+(t))$
is a sub-stochastic random walk on the half-space $\Z^{d-1}\times\N^*$ having transition matrix
\[
\left(p(z,z')=\mu(z'-z), \; z,z'\in\Z^{d-1}\times\N^*\right).
\] 
The random walk   $(Z_+(t))$ is identical to the homogeneous random walk  on
$\Z^d$ killed upon hitting the boundary  hyper-plane $\Z^{d-1}\times\{0\}$.

\subsection{Sample path large deviation principle for scaled processes} 
Before to formulate our large deviation results we recall the definition of the sample path large
deviation principle. 

\smallskip
\noindent
{\bf Definitions : \;}{\em 1) Let $D([0,T],\R^{d})$ denote the set of all right continuous with left
limits functions from $[0,T]$ to $\R^{d}$ endowed with Skorohod metric
(see Billingsley~\cite{Billingsley}). Recall that   
a mapping $I_{[0,T]}:~D([0,T],\R^{d})\to
[0,+\infty]$ is a good rate function on $D([0,T],\R^{d})$ if for any $c\geq 0$ and 
 any compact set $V\subset \R^{d}$, the set
\[
\{ \varphi \in D([0,T],\R^{d}): ~\phi(0)\in V \; \mbox{
and } \; I_{[0,T]}(\varphi) \leq c \}
\]
is compact in $D([0,T],\R^{d})$.  According to this definition, a good
rate function is lower semi-continuous. 

2) For a Markov chain $(Z(t))$ on $E\subset\R^d$ the family of scaled processes
$(Z^\eps(t) =\eps Z([t/\eps]), 
\,t\in[0,T])$,  is said to
satisfy {\it sample path large deviation principle} in $D([0,T], \R^{d})$ with a rate function
$I_{[0,T]}$ if for any $z\in\R^{d}$ 
\begin{equation}\label{e5-1}
\lim_{\delta\to 0} \;\liminf_{\eps\to 0} \; \inf_{z'\in E : |\eps z'-z|<\delta} \eps
\log\P_{z'}\left( Z^\eps(\cdot)\in {\cal 
O}\right) \geq -\inf_{\phi\in{\cal O}:\phi(0)=z} I_{[0,T]}(\phi), 
\end{equation}
for every open set ${\cal
O}\subset D([0,T],\R^{d})$,
and
\begin{equation}\label{e5-2}
\lim_{\delta\to 0} \;\limsup_{\eps\to 0} \; \sup_{z' \in E : |\eps z'-z|<\delta}
\eps\log\P_{z'}\left( Z^\eps(\cdot)\in 
F\right) \leq -\inf_{\phi\in F:\phi(0)=z} I_{[0,T]}(\phi).
\end{equation}
 for every closed set $F\subset   D([0,T],\R^{d})$. 
}

We refer to sample path large deviation
principle as SPLD principle. Inequalities (\ref{e5-1}) and (\ref{e5-2}) are
referred as lower and upper SPLD bounds respectively.

Proposition~4.1 of Ignatiouk~\cite{Ignatiouk:06} proves that under the
hypotheses (H2) and (H4), the family of
scaled processes $(Z_+^\eps(t)=\eps Z_+([t/\eps]), \, t\in[0,T])$ satisfies SPLD principle in
$D([0,T],\R^{d})$ with a 
good rate function 
\[
I_{[0,T]}^+(\phi) ~=~ \begin{cases} \int_0^T (\log \varphi)^*(\dot\phi(t)) \, dt, &\text{ if
    $\phi$ is absolutely continuous and }\\
&\text{\; $\phi(t)\in\R^{d-1}\times\R_+$ for all $t\in[0,T]$,}\\ 
+\infty &\text{ otherwise}
\end{cases}
\]
where $(\log \varphi)^*$ denotes the convex conjugate of the function $\log \varphi$ : 
\[
(\log \varphi)^*(v) ~\dot=~ \sup_{a\in\R^{d}} \Bigl(a\cdot v - \log
\varphi(a)\Bigr). 
\]
The next 
proposition provides the SPLD principle for the  scaled processes
$Z^\eps(t)$. 
\begin{prop}\label{pr5-1} Under  the hypotheses $(H_0)-(H_4)$, for every $T>0$,  the family of  scaled  
  processes $(Z^\eps([t/\eps]) ~=~ \eps Z([t/\eps]), \, t\in[0,T])$ satisfies SPLD
  principle in $D([0,T], \R^{d})$  with a 
  good rate function 
\[
I_{[0,T]}(\phi) ~=~ \begin{cases} \int_0^T L(\phi(t),\dot\phi(t)) \, dt, &\text{ if
    $\phi$ is absolutely continuous and  }\\ &\text{ $\phi(t)\in\R^{d-1}\times\R_+$ for all
    $t\in[0,T]$,}\\
+\infty &\text{ otherwise.}
\end{cases}
\]
The local rate function $L(z,v)$ is defined for every $z=(x,y), v\in\R^{d-1}\times\R$ by
the equality 
\[
L(z,v) ~=~ \begin{cases} (\log \varphi)^*(v) &\text{ if $y>0$}\\
(\log~\max\{\varphi, \varphi_0\})^*(v) &\text{ if $y=0$}
\end{cases}
\]
where $(\log\max\{\varphi,\varphi_0\})^*$ is
the convex conjugate of the function 
$\log \max\{\varphi, \varphi_0\}$~: 
\[
(\log \max\{\varphi, \varphi_0\})^*(v) ~=~ \sup_{a\in\R^{d}} ~\left( a\cdot v -
\log~\max\{\varphi(a), \varphi_0(a)\} \right).
\]
\end{prop}
\noindent

This proposition is a consequence of the results obtained in
\cite{D-E-W,D-E,Ignatiouk:02,Ignatiouk:04}. The results of Dupuis, Ellis and
Weiss~\cite{D-E-W} prove that $I_{[0,T]}$ is a good rate function on $D([0,T],\R^d)$
and provide the SPLD upper bound. Because of the communication condition, SPLD lower
bound follows from the local estimates 
obtained in \cite{Ignatiouk:02}, the general SPLD lower bound of Dupuis and
Ellis~\cite{D-E} and the integral representation of the corresponding rate function
obtained in \cite{Ignatiouk:04}. For the related results, see also \cite{D-B,
  D-E-2,I-M-S,S-W}. 

\medskip 

\subsection{Explicit form of quasi-potentials}
For a given rate function $J_{[0,T]}$ on  the Skorohod space $D([0,T],\R^d)$, the quantity 
\[
J(q,q') ~=~ \inf_{T>0} ~\inf_{\substack{\phi\in D([0,T],\R^d) : \\ \phi(0)=q, \phi(T) =
    q'}} J_{[0,T]}(\phi) 
\]
represents the optimal large deviation cost to go from $q$ to $q'$. Following Freidlin and
Wentzel terminology~\cite{W-F}, such a function $I : \R^d\times\R^d\to \R_+$ is called
{\em quasi-potential}. Borovkov 
and Mogulskii~\cite{Borovkov:01} called this function {\em second deviation rate
  function}. 

In this section, we calculate explicitly the quasi-potentials 
 $I(0,q)$ and $I^+(q',q)$  of the rate functions $I_{[0,T ]}$ and $I^+_{[0,T]}$ respectively.

\begin{prop}\label{pr5-2}  Under the hypotheses (H2) and (H4), for any $q',q\in \R^{d-1}\times\R_+$,
\be\label{e5-3} 
I^+(q',q) ~=~ \sup_{a\in D} a\cdot (q-q').
\ee
\end{prop}
\begin{proof} Indeed, 
for any $T>0$ and any
  absolutely continuous function $\phi :[0,T]\to\R^{d-1}\times\R_+$ with $\phi(0)=q'$
and $\phi(T)=q$,   
\[
I_{[0,T]}^+(\phi) ~=~ \int_0^T (\log \varphi)^*(\dot\phi(t)) \, dt ~\geq~ T (\log
\varphi)^*\left(\frac{q-q'}{T}\right)  
\]
because the function $(\log \varphi)^*$  is convex. Since the last relation holds
with the equality for the linear function $\phi(t)=t (q-q')/T$, $t\in[0,T]$ we  obtain   
\be\label{e5-4}
I^+(0,q) ~=~ \inf_{T>0} T (\log \varphi)^*\left(\frac{q-q'}{T}\right). 
\ee
Furthermore, under the hypotheses (H2) and (H4),
the function $\log \varphi$ is convex and continuous on $\R^d$ and hence,  it  is a closed
convex proper function on $\R^d$.   By Theorem~13.5  of
Rockafellar~\cite{R}  from this it follows that the support function of the set $D=\{a\in\R^d :~ \log\varphi(a)\leq
0\}$ 
is equal to the closure of the positively homogeneous convex
function $k$ generated by  $(\log \varphi)^*$. For any $v\in\R^d$ we have therefore  
\[
\text{cl} (k)(v) ~=~\sup_{a\in D} a\cdot v
\]
Moreover, under the hypotheses (H2) and (H4),
$(\log \varphi)^*$ is also a closed convex proper function on $\R^d$ with  
\[
0 ~<~ (\log \varphi)^*(0) ~\dot=~ - \inf_{a\in\R^d} \varphi(a) ~<~ +\infty. 
\]
By Theorem~9.7 of Rockafellar~\cite{R} from
this it follows that  the  positively homogeneous convex
function $k$ generated by  $(\log \varphi)^*$ is closed and for any
$q'q\in\R^{d-1}\times\{0\}$, the quantity $k(q-q')$ is equal to the right hand side
of \eqref{e5-4}. Hence, for any $q',q\in\R^{d-1}\times\{0\}$, 
\[
I^+(q',q) ~=~ k(q-q') ~=~ \text{cl} (k)(q-q') ~=~ \sup_{a\in D} a\cdot (q-q'). 
\]
Proposition~\ref{pr5-2} is therefore proved. 
\end{proof}

The next proposition identifies the quasi-potential of the rate function $I_{[0,T]}$. 

\begin{prop}\label{pr5-3}  Under the hypotheses (H0)-(H4),  for any non-zero vector 
 $q\in \R^{d-1}\times\R_+$, 
\be\label{e5-5}
I(0,q) ~=~ \inf_{\gamma\in\R^{d-1}\times\{0\}} I(0,\gamma) + I^+(\gamma,q)  ~=~
\sup_{a\in(\Theta\times\R)\cap D} ~a\cdot q  
\ee
\end{prop} 
\begin{proof} Indeed, the first equality of \eqref{e5-5} holds because for  any absolutely continuous function
  $\phi : [0,T] \to \R^{d-1}\times\R_+$ with  $\tau = \sup\{t > 0 :
  \phi(t)\in\R^{d-1}\times\{0\}\}$, one has  
\begin{align*}
I_{[0,T]}(\phi) &~=~ \int_0^T L(\phi(t),\dot\phi(t)) \, dt ~=~ \int_0^\tau 
L(\phi(t),\dot\phi(t)) \, dt + \int_\tau^T (\log \varphi)^*(\dot\phi(t)) \, dt\\
&~=~ I_{[0,\tau]}(\phi) + I^+_{[0,T-\tau]}(\phi_\tau) 
\end{align*}
where $\phi:[0,\tau]\to\R^{d-1}\times\R_+$ is the restriction of the function $\phi$ on
$[0,\tau]$ and $\phi_\tau:[0,T-\tau]\to\R^{d-1}\times\R_+$ is defined by
$\phi_\tau(t)=\phi(\tau+t)$ for all $t\in[0, T-\tau]$. 
To get the second equality of
\eqref{e5-5} we first  notice that  for any $T>0$ and any absolutely 
continuous function   $\phi:[0,T]\to\R^{d-1}\times\R_+$ with $\phi(0)=0$ 
and $\phi(T)= \gamma\in\R^{d-1}\times\{0\}$, the following relations hold   
\begin{align*}
I_{[0,T]}(\phi) &=~ \int_0^T L(\phi(t),\dot\phi(t)) \, dt \geq~ \int_0^T
(\log~\max\{\varphi, \varphi_0\})^*(\dot\phi(t)) \, dt \\ &\geq~  T (\log~\max\{
\varphi, \varphi_0\})^*\left(\frac{\gamma}{T}\right).
\end{align*}  
The first inequality holds here because according to the definition of the local rate
function, 
\[L(x,v) ~\geq~ (\log\max\{\varphi,
\varphi_0\})^*(v), \quad  \forall v\in\R^d, \; x\in\R^{d-1}\times\R_+.
\]
The second inequality 
is satisfied because the function $(\log~\max\{ \varphi, 
\varphi_0\})^*$ is convex. Since these relations hold with the equalities for the linear
function $\phi(t) = t \gamma/T$, we obtain 
\[
I(0,\gamma) ~=~ \inf_{T>0} T (\log~\max\{
\varphi, \varphi_0\})^*\left(\frac{\gamma}{T}\right), \quad \forall \gamma\in\R^{d-1}\times\{0\}, 
\]
and using next the same arguments as in the proof of Proposition~\ref{pr5-2} we conclude
that 
\[
I(0,\gamma) ~=~ \sup_{a :~ \max\{\varphi(a),\varphi_0(a)\}\leq 1} a\cdot \gamma ~=~ \sup_{a\in D\cap
  D_0} a\cdot \gamma, \quad \quad \forall \gamma\in\R^{d-1}\times\{0\}. 
\]
From the last relation it follows that 
\[
I(0,\gamma) ~=~ \sup_{a\in\Theta\times\R} a\cdot\gamma, \quad \quad \forall \gamma\in\R^{d-1}\times\{0\} 
\]
because  the set $\Theta\times\{0\}$ is the orthogonal projection of the set $D\cap D_0$ onto
the hyper-plane $\R^{d-1}\times\{0\}$. Using Proposition~\ref{pr5-3} we obtain therefore 
\[
\inf_{\gamma\in\R^{d-1}\times\{0\}} I(0,\gamma) + I^+(\gamma, q) ~=~
\inf_{\gamma\in\R^{d-1}\times\{0\}} \left(\sup_{a\in\Theta\times\R} a\cdot\gamma +
\sup_{a\in D} a\cdot(q-\gamma)\right) 
\]
Moreover, since for $\gamma\in\R^d$ with a non-zero last coordinate on has 
\[
\sup_{a\in\Theta\times\R} a\cdot\gamma ~=~ + \infty,
\]
the infimum over $\gamma\in\R^{d-1}\times\{0\}$ at the right hand side of the above
relation can be replaced by the infimum over $\gamma\in\R^d$. Finally, under the hypotheses
(H0)-(H4), the interior of the set $\hat{D}=(\Theta\times\R)\cap D$ is non-empty and
consequently, by Corollary~16.4.1 of Rockafellar~\cite{R}, 
\[
\inf_{\gamma\in\R^d} \left(\sup_{a\in\Theta\times\R} a\cdot\gamma +
\sup_{a\in D} a\cdot(q-\gamma)\right)  ~=~ \sup_{a\in (\Theta\times\R)\cap D} a\cdot q. 
\]
The second equality of \eqref{e5-5} is therefore proved. 
\end{proof}

\begin{cor}\label{cor5-1} Under the hypotheses (H0) - (H4), the functions 
  $q\to I^+(0,q)$ and $q\to I(0,q)$ are convex and continuous everywhere on $\R^{d-1}\times\R_+$. 
\end{cor}
\begin{proof} Indeed, the equalities \eqref{e5-3} and \eqref{e5-5} show that each of these
  functions is a support function of a compact set. From this it follows that they 
  are finite, convex and therefore continuous on $\R^{d-1}\times\R_+$.  
\end{proof}

The next proposition investigates the point where the minimum of the function $\gamma \to
I(0,\gamma) + I^+(\gamma,q)$  over $\gamma\in\R^{d-1}\times\{0\}$ is attained. 
Recall that by Corollary~\ref{lem2-5}, for every
$q\in\R^{d-1}\times\R_+^*$, there exists a unique point $\gamma_q\in\R^{d-1}\times\{0\}$
for which the vectors $\gamma_q, \, q-\gamma_q$ are normal to the set
$\hat{D}=(\Theta\times\R)\cap D$  at the point $\hat{a}(q)$ and the vector $q-\gamma_q$ is
normal to the set $D$ at the point $\hat{a}(q)$.

\begin{prop}\label{pr5-4} For  $q\in\R^{d-1}\times\R_+^*$, 
the point $\gamma_q$ is the only 
minimum of the function $\gamma \to I(0,\gamma) + I^+(\gamma,q)$ on the hyperplane 
$\R^{d-1}\times\{0\}$.  
\end{prop}
\begin{proof} 
Indeed, by Corollary~\ref{cor5-1}, the functions $\gamma \to I(0,\gamma)$
and $\gamma\to I^+(\gamma, q) = I^+(0, q-\gamma)$ are  finite and convex everywhere on $\R^{d-1}\times\{0\}$. The function $\gamma
\to    I(0,\gamma) + I^+(\gamma,q)$ is therefore also finite and convex everywhere on
$\R^{d-1}\times\{0\}$. By Theorem~23.5 of
Rockafellar~\cite{R} from this it follows  that $I(0,\gamma) + I^+(\gamma, q)$
achieves its minimum over $\R^{d-1}\times\{0\}$ at the point $\hat{\gamma}\in\R^{d-1}\times\{0\}$ if and only if the
differential $\partial ( I(0,\hat\gamma) + I^+(\hat\gamma, q))$ of this function at the point 
$\hat{\gamma}$ contains zero vector. 
Theorem~23.8 of Rockafellar~\cite{R} proves that 
\[
\partial ( I(0,\hat\gamma) + I^+(\hat\gamma,q)) ~=~ \partial I(0,\hat\gamma) + \partial I^+(\hat\gamma,
q)
\]
where $\partial I(0,\hat\gamma)$ denotes the differential of the function
$\gamma\to I(0,\gamma)$ and $\partial I^+(\hat\gamma,
q)$ is the differential of the function
$\gamma\to I^+(\gamma,q)$ at the point $\gamma=\hat\gamma$. 
By Corollary~23.5.3 of Rockafellar~\cite{R}, from 
\eqref{e5-5} it follows that  $\hat{a}\in \partial I(0,\gamma)$ if and only
if $\hat{a}\in D$ and 
\[
I(0,\gamma) ~=~ \sup_{a\in(\Theta\times\R)\cap D} a\cdot \gamma ~=~ \hat{a}\cdot \gamma
\]
or equivalently, when the vector $\gamma$ is normal to
the set $(\Theta\times\R)\cap D$ at the point $\hat{a}\in (\Theta\times\R)\cap
D$. Similarly, from
\eqref{e5-3} it follows that  $a'= -\hat{a}\in\partial I^+(\gamma,q)$  
if and only if $\hat{a}\in D$ and 
\[
I^+(\gamma,q) ~=~ \sup_{a\in(\Theta\times\R)\cap D} a\cdot (q-\gamma) ~=~ \hat{a}\cdot (q-\gamma)
\]
or equivalently, when the vector $(q-\gamma)$ is normal to the set $D$ at the
point $\hat{a}\in D$.  According to the definition of $\gamma_q$, this proves that the
function $\gamma\to I(0,\gamma) + I^+(\gamma, q)$ achieves its minimum  over the set
$\R^{d-1}\times\{0\}$ at the point $\gamma_q\in\R^{d-1}\times\{0\}$.

Conversely, if  the
function $\gamma\to I(0,\gamma) + I^+(\gamma, q)$ achieves its minimum  over the set
$\R^{d-1}\times\{0\}$ at some point $\hat\gamma\in\R^{d-1}\times\{0\}$ then there is a
point $\hat{a}\in(\Theta\times\R)\cap D$ for which the following conditions are satisfied~:
\begin{itemize}
\item[--] the vector $\hat\gamma$ is normal to 
the set $(\Theta\times\R)\cap D$ at the point $\hat{a}$,
\item[--] the vector $(q-\hat\gamma)$ is normal to the set $D$ at the
point $\hat{a}$,
\item[--] and  $
I(0,\hat\gamma) + I^+(\hat\gamma, q) ~=~ \hat{a}\cdot \hat\gamma + \hat{a}(q-\gamma) ~=~
\hat{a}\cdot q$.
\end{itemize}
Moreover, from 
\eqref{e5-5} it follows that 
\[
I(0,\hat\gamma) + I^+(\hat\gamma, q) ~=~ \sup_{a\in(\Theta\times\R)\cap D} a\cdot q 
\]
and consequently, 
\[
\sup_{a\in(\Theta\times\R)\cap D} a\cdot q ~=~ \hat{a}\cdot q
\]
The last relation shows that the vector $q$ is normal to the set $(\Theta\times\R)\cap D$ at the point
$\hat{a}\in(\Theta\times\R)\cap D$. By Lemma~\ref{lem2-3}  from this it follows that
$\hat{a}= \hat{a}(q)$. The vectors $\gamma_q, \, q-\gamma_q$ are therefore normal to the set
$\hat{D}=(\Theta\times\R)\cap D$  at the point $\hat{a}(q)$ and the vector $q-\gamma_q$ is
normal to the set $D$ at the point $\hat{a}(q)$. By Corollary~\ref{lem2-5} this proves that $\hat{\gamma}=\gamma_q$. 
\end{proof}

\subsection{Logarithmic asymptotics of Green's function}\label{sec-laGf} Now, we
  obtain logarithmic asymptotics of  Green's functions $G(z,z')$ and $G_+(z,z')$ for
  the Markov processes $(Z(t))$ and $(Z_+(t))$. 

\begin{prop}\label{pr5-5} Under the hypotheses (H2)-(H4), for any
  $q\in\R^{d-1}\times\R_+$ 
  and any sequences $\eps_n >0$ and  
  $z_n\in\Z^{d-1}\times\N^*$ with $\lim_n \eps_n = 0$ and $\lim_n \, \eps_n z_n = q$  the following   
  relations hold 
\[
\lim_{n\to\infty} ~\eps_n\log G_+(z,z_n)  ~=~ - \sup_{a\in D} a\cdot q.
\]
\end{prop}
\begin{proof} Indeed, let the sequences $\eps_n >0$ and  $z_n\in\Z^{d-1}\times\N^*$ be
  such that $\lim_n \eps_n =0$ and 
  $\lim_n\eps_n z_n = q$. Then by Lemma~\ref{lem2-1}, for  any $a\in D$,  
\[
G_+(z,z_n) ~\leq~ \exp(a\cdot(z-z_n)) G_S(0,0), \quad \quad \forall n\in\N,
z\in\Z^{d-1}\times\N^*
\]
and consequently, 
\[
\lim_{n\to\infty} ~\eps_n\log G_+(z,z_n)  ~\leq~ - a\cdot q \quad \quad \forall a\in D
\]
from which it follows that
\[
\lim_{n\to\infty} ~\eps_n\log G_+(z,z_n)  ~\leq~ - \sup_{a\in D}a\cdot q. 
\]
The  inequality 
\[
\lim_{n\to\infty} ~\eps_n\log G_+(z,z_n)  ~\geq~ - \sup_{a\in D}a\cdot q 
\]
was proved in Proposition~4.2 of 
Ignatiouk~\cite{Ignatiouk:06} by using lower large deviation bound  for the 
scaled processes $Z_+^{\eps}(t) = \eps Z_+(t/\eps)$ and communication condition. This
proof is quite similar to the proof of the lower bound \eqref{e5-6} below.\end{proof}  

\begin{prop}\label{pr5-6} Under the hypotheses (H0)-(H4), for any
  $q\in\R^{d-1}\times\R_+$, 
  and any sequences $\eps_n >0$ and  
  $z_n\in\Z^{d-1}\times\N$ with $\lim_n \eps_n = 0$,  and $\lim_{n\to\infty}\eps_n z_n= q$  the following   
  relation holds : 
\[
\lim_{n\to\infty} ~\eps_n \log G(z,z_n)    ~=~ - \sup_{a\in(\Theta\times\R)\cap D} a \cdot q, \quad \quad \forall  z\in\Z^{d-1}\times \N.
\]
\end{prop}
\begin{proof}  Let two sequences $\eps_n>0$ and $z_n\in\Z^{d-1}\times\N$ be such that
  $\lim_n\eps_n = 0$ and $\lim_n \eps_n z_n = q$. We begin our analysis with the proof of the lower bound 
\be\label{e5-6}
\lim_{n\to\infty} ~\eps_n \log G(z,z_n) ~\geq~ - \sup_{a\in(\Theta\times\R)\cap D} a \cdot q. 
\ee 
For this we use the lower large deviation bound and communication condition. Denote for
$B\in\R^d$   
\[
G(z, B) ~=~ \sum_{z'\in B\cap\Z^{d-1}\times\N} G(z,z'). 
\] 
The large deviation lower bound implies that for any $\delta > 0$ and $T>0$, 
\begin{align*}
  ~\liminf_{n\to\infty} ~\eps\log  G(z,\eps^{-1} B(q,\delta))  &~\geq~
  ~\liminf_{n\to\infty} ~\eps\log  \P_z\left(Z_\eps(T) \in B(q,\delta)\right)   \\ &~\geq~
 -  \inf_{\substack{\phi\in D([0,T],\R^{d-1}\times\R_+) :~ \phi(0) = 0,
     \; \phi(T) \in B(q,\delta)} } I_{[0,T]}(\phi) 
\\ &~\geq~ - \inf_{\substack{\phi\in D([0,T],\R^{d-1}\times\R_+) :~ \phi(0) = 0, 
     \; \phi(T) =q } } I_{[0,T]}(\phi)  
\end{align*} 
from which it follows that  
\be\label{e5-7}
 \lim_{\delta\to 0}
  ~\liminf_{n\to\infty} ~\eps\log  G(z,\eps^{-1} B(q,\delta))   ~\geq~ - I(0,q) =
 - \sup_{a\in(\Theta\times\R)\cap D} a \cdot q
\ee
where the last relation  is proved by 
Proposition~\ref{pr5-3}. Moreover, by
Proposition~\ref{pr4-1}, the Markov process $(Z(t))$ satisfies communication 
condition and hence, there are $\theta>0$ and $C>0$ such that for any 
$z',z''\in\Z^{d-1}\times\N$ such that $z'\not=z''$, the probability that the Markov process
$(Z(t))$ starting at $z'$  hits $z''$ before the first return to $z'$ is greater than
$\theta^{C|z''-z'|}$. This proves that for any
  $z,z'\in\Z^{d-1}\times\N$ and  $n\in\N$ 
\[
G(z,z_n) ~\geq~ G(z,z')
\theta^{C|z_n-z'|} 
\]
and consequently, for all those $n\in\N$ for which $|q- \eps_n z_n| <\delta$, we obtain 
\begin{align*}
G(z, \eps_n^{-1}B(q,\delta)) ~\theta^{2C\delta /\eps_n} &~\leq~ \sum_{z'\in \Z^{d-1}\times\N :~
  z'\in \eps_n^{-1}B(q,\delta)} G(z,z')
\theta^{C|z_n-z'|}  \\&~\leq~ \text{Card}\{z\in\Z^d : z\in
\eps_n^{-1}B(q,\delta) \} ~G(z,z_n) \\ &~\leq~ (2\delta \eps_n^{-1} +1)^d ~G(z,z_n)
\end{align*}
The last inequality shows that 
\[
\lim_{n\to\infty} ~\eps_n\log G(z,z_n)  ~\geq~ 2C \delta \log \theta + \liminf_{\eps\to 0} ~\eps
~\log
G\bigl(z, \eps^{-1}B(q,\delta) \bigr)
\]
and hence, letting $\delta\to 0$ and using \eqref{e5-7}, we get \eqref{e5-6}

To prove the inequality 
\be\label{e5-8}
\lim_{n\to\infty} ~\eps_n \log G(z,z_n)   ~\leq~ - \sup_{a\in(\Theta\times\R)\cap D} a \cdot q
\ee
 we use Lemma~\ref{lem2-1} and Corollary~\ref{cor2-2}. 
For $a\in D$, $z\in\Z^{d-1}\times\N$ and
$z'\in\Z^{d-1}\times\{0\}$, by Lemma~\ref{lem2-1},   
\[
G(z,z') ~\leq~  G(z',z') \exp(\ol{a} \cdot z -
\ol{a}\cdot z') ~=~ G(0,0) \exp(\ol{a} \cdot z -
a\cdot z') . 
\]
Moreover, if $\varphi_0(\ol{a})< 1$ then for $z\in\Z^{d-1}\times\N$ and
$z'\in\Z^{d-1}\times\N^*$, by Corollary~\ref{cor2-2},  
\[
\frac{G(z,z')}{G_S(0,0)} \leq \exp(a\cdot (z-z')) + \varphi_0(a)
(1-\varphi_0(\ol{a}))^{-1} \exp(\ol{a} \cdot z - a\cdot z').
\]
These  inequalities  show that for any $a\in D$ for which $\varphi_0(\ol{a}) < 1$, one has 
\[
\lim_{n\to\infty} ~\eps_n \log G(z,z_n)   ~\leq~ -   
a\cdot q.
\]
The last relation proves \eqref{e5-8} because  $(\Theta\times\R)\cap D = \{a\in D
: \varphi_0(\ol{a}) \leq 1\}$. \end{proof}

\section{Principal part of the renewal equation}\label{sec-reneq}   Recall that  the transition probabilities
 $p(z,z')$ of the Markov process $(Z(t))$
 are the same as transition probabilities $p(z,z')=\mu(z'-z)$  
 of the homogeneous random walk $S(t)$  on $\Z^d$  for $z\in\Z^{d-1}\times\N_+^*$ and that 
 $p(z,z')=\mu_0(z'-z)$ for $z\in\Z^{d-1}\times\{0\}$.  From this it follows that
 the Green's function $G(z,z')$ 
 satisfies  the following renewal equation 
\be\label{e6-1}
G(z,z_n) ~=~ G_+(z,z_n) ~+ \sum_{\substack{w\in
    \Z^{d-1}\times\{0\}\\ w'\in\Z^{d-1}\times\N^*}} G(z,w) \mu_0(w'-w) G_+(w',z_n). 
\ee
$G_+(z,z')$ denotes here Green's function  of the homogeneous random
 walk $(Z_+(t))$ killed upon 
 hitting the half-space $\Z^{d-1}\times(-\N)$ :  
\[
G_+(z,z') ~\dot=~ \sum_{t=0}^\infty \P_z(Z_+(t)=z') ~\dot=~ \sum_{t=0}^\infty
\P_z(S(t)=z'; \, \tau > t)
\]
where $\tau = \inf\{t\geq 1 : S(t)\in\Z^{d-1}\times\{0\}\}$.  

In this section we show that for a sequence $z_n\in\Z^{d-1}\times\N$
with $\lim_n |z_n|=\infty$ and $\lim_n z_n/|z_n|= q\in\R^{d-1}\times]0,+\infty[$  the right hand side of the renewal equation
\eqref{e6-1} can be decomposed into a main part 
\[
\Xi_{\delta}^q(z,z_n) ~=~ G_+(z,z_n) ~\1_{\{\gamma_q=0\}} + \hspace{-4mm}\sum_{\substack{w\in
    \Z^{d-1}\times\{0\}, w'\in \Z^{d-1}\times\N^* :\\ |w-\gamma_q |z_n||\leq \delta|z_n|
    }}\hspace{-8mm} G(z,w)  \mu_0(w'-w) G_+(w',z_n)
\]
and the corresponding negligible part $G(z,z_n) - \Xi_{\delta}^q(z,z_n)$.  Recall that $\gamma_q$ is a only vector on the
   boundary hyperplane $\R^{d-1}\times\{0\}$ for which the vectors 
   $\gamma_q$ and $q-\gamma_q$ belong to the normal cone $V(\hat{a}(q))$ to the set
   $\hat{D} ~\dot=~\{a\in D : \varphi'\ol{a}) \leq 1\} = (\Theta\times\R)\cap D$ at the point  $\hat{a}(q)$  and the vector $q-\gamma_q$ is
   normal to the set $D$ at the point $\hat{a}(q)$ (see Corollary~\ref{cor2-4}). By
   Propositions~\ref{pr5-3} and ~\ref{pr5-4}, this is also
   the only minimum of the function $\gamma\to I(0,\gamma) + I^+(\gamma,q)$ on the
   boundary hyperplane $\R^{d-1}\times\{0\}$ where 
\[
I(0,q) ~=~ I(0,\gamma_q) +
   I^+(\gamma_q,q).
\]
We begin our analysis we the following lemma. 

\begin{lemma}\label{lem6-1} Under the hypotheses (H0)-(H4), 
\[
I_{min} ~\dot=~ \inf_{\gamma\in\R^{d-1}\times\{0\}:~ |\gamma|=1} I(0,\gamma) + I^+(\gamma,0) ~>~ 0
\]
\end{lemma}
\begin{proof} Indeed, by Corollary~\ref{cor5-1}, the function 
\[
\gamma\to I(0,\gamma) +
  I^+(\gamma,0) = I(0,\gamma) + I^+(0,-\gamma)
\]
is continuous. To prove our lemma it is therefore sufficient to show that  
\be\label{e6-2}
I(0,\gamma) + I^+(0,-\gamma) > 0, \quad \quad \forall \; \gamma\not=0.
\ee
To prove this inequality let us notice that 
\begin{align*}
I(0,\gamma) + I^+(0,-\gamma) &~=~ \sup_{a\in\hat{D}} a\cdot\gamma + \sup_{a\in
  D}a\cdot(-\gamma) \\&~\geq~ \sup_{a\in\hat{D}} a\cdot\gamma + \sup_{a\in
  \hat{D}}a\cdot(-\gamma)  ~=~ \sup_{a\in\hat{D}} a\cdot\gamma ~-~ \inf_{a\in
  \hat{D}}a\cdot\gamma ~\geq~ 0
\end{align*}
where the last relation holds with equality if and only if $a\cdot\gamma=0$ for all
$a\in\hat{D}$. Under the hypotheses of our lemma, for any non-zero vector $\gamma\in\R^d$
there is $a\in\hat{D}$ for which $a\cdot\gamma \not= 0$ because the set $\hat{D}$ has a non-empty
interior (see the proof of Lemma~\ref{lem2-5}). The inequality \eqref{e6-2} is therefore proved. 
\end{proof}

\begin{prop}\label{pr6-1} Under the hypotheses (H0)-(H4), for any
   $q\in{\cal S}^d_+\cap \R^{d-1}\times]0,+\infty[$ $\delta >0$  and any sequence 
  $z_n\in\Z^{d-1}\times\N$ with $\lim_n|z_n|=\infty$ and $\lim_nz_n/|z_n|=q$, 
\be\label{e6-3}
\lim_{n\to\infty} ~\Xi_{\delta}^q(z,z_n)/G(z,z_n) ~=~ 1, \quad \quad \forall  z\in\Z^{d-1}\times\N.
\ee  
\end{prop}
\begin{proof} Let $q\in{\cal S}^d_+\cap \R^{d-1}\times]0,+\infty[$ and let a sequence  $z_n\in
  \Z^{d-1}\times\N$ be such that $|z_n|\to\infty$ and 
  $z_n/|z_n|\to q$ as $n\to\infty$. Then by Proposition~\ref{pr5-4}, 
\[
\lim_{n\to\infty} \frac{1}{|z_n|} \log G(z,z_n) ~=~ - ~I(0,q),
\]
and hence, to get \eqref{e6-3} it is sufficient to show that 
\[
\limsup_{n\to\infty} \frac{1}{|z_n|} \log(G(z,z_n) - \Xi_{\delta}^q(z,z_n)) ~<~ - ~I(0,q). 
\]
By  Lemma~1.2.15 of Dembo and Zeitouni~\cite{D-Z}, for this it is sufficient to prove the
following three inequalities~:
\be\label{e6-4}
\limsup_{n\to\infty} \frac{1}{|z_n|} \log G_+(z,z_n)  < - I(0,q) \quad \text{ when } \quad
\gamma_q \not= 0,
\ee
\be\label{e6-5}
\limsup_{n\to\infty} \frac{1}{|z_n|} ~\log \hspace{-4mm} \sum_{\substack{w\in
    \Z^{d-1}\times\{0\}, w'\in \Z^{d-1}\times\N^* :\\  
    ~|w'-w| > \delta' |z_n|
    }}\hspace{-7mm} G(z,w)  \mu_0(w'-w) G_+(w',z_n) ~<~ - I(0,q) 
\ee 
and 
\be\label{e6-6}
\limsup_{n\to\infty} \frac{1}{|z_n|} ~\log \hspace{-4mm} \sum_{\substack{w\in
    \Z^{d-1}\times\{0\}, w'\in \Z^{d-1}\times\N^* :\\ 
    |w-\gamma_q |z_n|| > \delta|z_n|, ~|w'-w| \leq \delta' |z_n|
    }}\hspace{-7mm} G(z,w)  \mu_0(w'-w) G_+(w',z_n) ~<~ - I(0,q)
\ee
for some $\delta' > 0$ small enough. 

{\em Proof of ~\eqref{e6-4}:}
This relation is a consequence of  Propositions~\ref{pr5-3},~\ref{pr5-4} and
~\ref{pr5-5}. Namely, Proposition~\ref{pr5-5} proves that 
\[
\lim_{n\to\infty} \frac{1}{|z_n|} \log G_+(z,z_n) ~=~ - ~I^+(0,q)   
\]
and by Propositions~\ref{pr5-3} and ~\ref{pr5-4},  
\[
I^+(0,q) ~=~ I(0,0) + I^+(0,q) ~>~
I(0,\gamma_q) + I^+(\gamma_q,q) ~=~ I(0,q) \quad \text{when } \quad \gamma_q\not=0.
\]
 
{\em Proof of ~\eqref{e6-5}:} For any
$a\in D$ for which $\varphi_0(a) < 1$, with the same arguments as in the proof of
Corollary~\ref{cor2-2} one gets  
\begin{align*}
\sum_{\substack{w\in
    \Z^{d-1}\times\{0\}, w'\in \Z^{d-1}\times\N^* :\\  
    ~|w'-w| > \delta' |z_n|
    }}&\hspace{-7mm} G(z,w)  \mu_0(w'-w) G_+(w',z_n)  \exp(a\cdot z_n) \\&~\leq~
G_S(0,0) \hspace{-4mm} \sum_{\substack{w\in
    \Z^{d-1}\times\{0\}, w'\in \Z^{d-1}\times\N^* :\\  
    ~|w'-w| > \delta' |z_n|
    }}\hspace{-7mm} G(z,w)  \mu_0(w'-w) \exp(a\cdot w')  \\&~\leq~ \frac{G_S(0,0)
\exp(\ol{a}\cdot z)}{1-\varphi_0(\ol{a})} \sum_{u\in \Z^{d-1}\times\N^* : 
    ~|u| > \delta' |z_n|
    }\hspace{-7mm}  \mu_0(u) \exp(a\cdot u). 
\end{align*}
Hence, the right hand side of \eqref{e6-5} does not exceed 
\[
- \lim_{n\to\infty} a\cdot z_n/|z_n| + \delta' \limsup_{R\to\infty} ~\frac{1}{R}\log \!\sum_{u \in\Z^{d-1}\times\N : ~|u| > R} \mu_0(u)
\exp(a\cdot u ) 
\]
where 
\[
\limsup_{R\to\infty} ~\frac{1}{R}\log \!\sum_{z \in\Z^{d-1}\times\N : ~|z| > R} \mu_0(z)
\exp(a\cdot z ) ~=~ - \infty
\]
because  under the hypotheses (H4), the function $a'\to \varphi_0(a+a')$ is finite
everywhere in $\R^d$.  Relation \eqref{e6-5} is therefore proved. 

{\em Proof of ~\eqref{e6-6}:}  
Lemma~\ref{lem2-1} proves that  for any 
$w\in\Z^{d-1}\times\{0\}$, $w'\in\Z^{d-1}\times\N^*$ and $a,a'\in
D$ with $\varphi_0(\ol{a}) \leq 1$ and $|w'-w|\leq \delta'|z_n|$, 
\[
G(z,w)G_+(w',z_n) ~\leq~ G(w,w) G_S(0,0) \exp(\ol{a}\cdot(z-w) +  a'\cdot(w'-z_n)) 
\]
where $G(w,w) = G(0,0)$ and 
\begin{align*}
\ol{a}\cdot z +  a'\cdot(w'-z_n) &= \ol{a}\cdot z +  a'\cdot(w- q|z_n|) + a'(q|z_n| - z_n) +
a'(w'-w) \\&\leq c|z| + c \bigl|q|z_n| - z _n\bigr| + \delta' c |z_n|
\end{align*}
with $c=\max_{a\in D} |a|$. Moreover, according to the definition of the mapping $a\to \ol{a}$, 
\[
\ol{a}\cdot w = a\cdot w
\]
because $w\in\Z^{d-1}\times\{0\}$ and consequently,    
\begin{align*}
G(z,w)G_+(w',z_n) &\leq~ G(0,0) G_S(0,0) \exp( - a\cdot w -
a'\cdot(q|z_n| - w)) \\ & \hspace{4cm} \times\exp( c |z_n - q|z_n||  + c |z| + 
\delta' c |z_n|)  
\end{align*}
Since the last inequality holds for arbitrary $a'\in D$ and $a\in\hat{D} ~\dot=~ \{a\in D :
\varphi(\ol{a})\leq 1\}$, 
using Propositions~\ref{pr5-2} and ~\ref{pr5-3}  we get 
\begin{align*}
G(z,w)G_+(w',z_n) &\leq  G(0,0) G_S(0,0) \exp(- I(0,w) - I^+(w,q|z_n|)) \\ & \hspace{4cm} \times\exp( c |z_n - q|z_n||  + c |z| +
\delta' c |z_n|)
\end{align*}
from which it follows that 
\begin{align*}
&\sum_{\substack{w\in \Z^{d-1}\times\{0\}, w'\in \Z^{d-1}\times\N^* :\\ 
   \delta|z_n| < |w-\gamma_q |z_n|| \leq R|z_n|, ~|w'-w| \leq \delta' |z_n|
    }}\hspace{-7mm} G(z,w)  \mu_0(w'-w) G_+(w',z_n)  \\ &\hspace{2.5cm}\leq\sum_{\substack{w\in
      \Z^{d-1}\times\{0\} :\\ 
   \delta|z_n| < |w-\gamma_q |z_n|| \leq R|z_n|, 
    }}\hspace{-7mm}  G(0,0) G_S(0,0) \exp(- I(0,w) - I^+(w,q|z_n|)) \\ &  \hspace{7cm} \times\exp( c |z_n - q|z_n||  + c |z| +
\delta' c |z_n|)  
\end{align*} 
and consequently, 
\begin{align*}
\limsup_{n\to\infty} \frac{1}{|z_n|} \log \hspace{-6mm} &\sum_{\substack{w\in
    \Z^{d-1}\times\{0\}, w'\in \Z^{d-1}\times\N^* :\\ 
   \delta|z_n| < |w-\gamma_q |z_n|| \leq R|z_n|, ~|w'-w| \leq \delta |z_n|
    }}\hspace{-7mm} G(z,w)  \mu_0(w'-w) G_+(w',z_n) \\ 
&\leq~ \limsup_{\eps\to 0} ~\eps \log  \hspace{-0.5cm}\sum_{\substack{w\in
      \Z^{d-1}\times\{0\} : 
   ~\delta < |\eps w-\gamma_q | \leq R , 
    }}\hspace{-7mm}  \exp(- I(0,w) - I^+(w,q/\eps)) + c\delta' \\ 
&\leq~ - \inf_{\gamma\in\R^{d-1}\times\{0\}: \delta < |\gamma_q-\gamma| \leq R} (I(0,\gamma) +
  I^+(\gamma,q) ) +  c\delta' 
\end{align*}
where the last relation holds because the number of points $w\in\Z^{d-1}\times\{0\}$
satisfying the inequality  $\delta < |\eps w-\gamma_q | \leq R$ does not exceed $(1 +
2R/\eps)^d$ and for each of them, 
$I(0,w) + I^+(w,q/\eps) ~=~ \eps^{-1} (I(0,\gamma) +
I^+(\gamma,q))$ with $\gamma = \eps w$. Recall now that the function $\gamma\to I(0,\gamma) +
I^+(\gamma,q)$ is convex and continuous on $\R^{d-1}\times\{0\}$, the point $\gamma_q$ is
the only minimum of this function at $\R^{d-1}\times\{0\}$ and  $
I(0,\gamma_q) + I^+(\gamma_q) = I(0,q)$ 
(see Corollary~\ref{cor5-1} and Propositions~\ref{pr5-3} and ~\ref{pr5-4}). This proves
that 
\begin{align*}
\inf_{\substack{\gamma\in\R^{d-1}\times\{0\}:\\ \delta < |\gamma_q-\gamma| \leq R}} (I(0,\gamma) +
  I^+(\gamma,q) ) &~\geq~ \inf_{\gamma\in\R^{d-1}\times\{0\}: ~\delta < |\gamma_q-\gamma| } (I(0,\gamma) +
  I^+(\gamma,q) )  \\ &~>~ I(0,\gamma_q) + I^+(\gamma_q,q) ~=~ I(0,q).  
\end{align*}
and consequently, for any $R>\delta > 0$ and  $\delta' >0$ satisfying the inequality  
\[
0 < \delta' c < \inf_{\gamma\in\R^{d-1}\times\{0\}: ~\delta < |\gamma_q-\gamma| } (I(0,\gamma) +
  I^+(\gamma,q) )  - I(0,q) 
\]
we get 
\[
\limsup_{n\to\infty} \frac{1}{|z_n|} \log \hspace{-6mm}  \sum_{\substack{w\in
    \Z^{d-1}\times\{0\}, w'\in \Z^{d-1}\times\N^* :\\ 
   \delta|z_n| < |w-\gamma_q |z_n|| \leq R|z_n|, ~|w'-w| \leq \delta' |z_n|
    }}\hspace{-7mm} G(z,w)  \mu_0(w'-w) G_+(w',z_n) ~<~ - I(0,q).
\]
Now, to complete the proof of \eqref{e6-6} it is sufficient to show that there is $R >
0$  such that 
\be\label{e6-7}
\limsup_{n\to\infty} \frac{1}{|z_n|} \log \hspace{-6mm} \sum_{\substack{w\in
    \Z^{d-1}\times\{0\}, w'\in \Z^{d-1}\times\N^* :\\ 
    |w-\gamma_q |z_n|| > R|z_n|, ~|w'-w| \leq \delta' |z_n|
    }}\hspace{-7mm} G(z,w)  \mu_0(w'-w) G_+(w',z_n) ~<~ - I(0,q).
\ee
To get this inequality we use again Lemma~\ref{lem2-1} combined with
Propositions~\ref{pr5-2} and ~\ref{pr5-3} :~ for any $a,a'\in D$, 
$w\in\Z^{d-1}\times\{0\}$ and $w'\in\Z^{d-1}\times\N^*$ with $\varphi(\ol{a})\leq 1$,
    and $|w-w'|\leq \delta'|z_n|$, from Lemma~\ref{lem2-1} it follows that   
\begin{align*}
G(z,w) G_+(w',z_n) &\leq G(w,w) G_S(0,0) \exp(\ol{a}\cdot (z-w) + a'\cdot(w'-z_n)) \\
&= G(w,w) G_S(0,0) \exp(- a\!\cdot \!w + a'\!\cdot \!w + \ol{a}\cdot z + a'\!\cdot \!(w'-w -z_n)) \\ 
&\leq G(w,w) G_S(0,0) \exp(- a\cdot w + a'\cdot w + c |z|  + (1+\delta') c
|z_n|) 
\end{align*}
with $c = \max_{a\in D}|a|$ and $G(w,w) = G(0,0)$. 
Using therefore Propositions~\ref{pr5-2} and ~\ref{pr5-3} we obtain 
\begin{align*}
G(z,w) G_+(w',z_n) &~\leq~G(0,0) G_S(0,0) \exp(- I(0,w) -  I^+(0,-w)) \\
&\hspace{5.5cm} \times \exp( c|z| + c (1+\delta') 
|z_n|)  
\end{align*}
from which it follows that 
\begin{align*}
&\hspace{-0.5cm}\sum_{\substack{w\in
    \Z^{d-1}\times\{0\}, w'\in \Z^{d-1}\times\N^* :\\ 
    |w-\gamma_q |z_n|| > R|z_n|, ~|w'-w| \leq \delta' |z_n|
    }}\hspace{-7mm} G(z,w)  \mu_0(w'-w) G_+(w',z_n) \\&~\leq~ \sum_{\substack{w\in
    \Z^{d-1}\times\{0\} :~
    |w-\gamma_q |z_n|| > R|z_n|, 
    }}\hspace{-7mm} G(0,0) G_S(0,0) \exp(- I(0,w) -  I^+(0,-w)) \\& \hspace{8cm} \times\exp(c|z| + c (1+\delta')
|z_n|)    
\end{align*}
and consequently, 
\begin{align*}
&\hspace{-0.5cm}\limsup_{n\to\infty} \frac{1}{|z_n|} \log \hspace{-4mm} \sum_{\substack{w\in
    \Z^{d-1}\times\{0\}, w'\in \Z^{d-1}\times\N^* :\\ 
    |w-\gamma_q |z_n|| > R|z_n|, ~|w'-w| \leq \delta' |z_n|
    }}\hspace{-7mm} G(z,w)  \mu_0(w'-w) G_+(w',z_n) \\ 
&\leq~ \limsup_{\eps\to 0} ~\eps
  \log \hspace{-5mm} \sum_{\gamma \in
    \eps \Z^{d-1}\times\{0\} :~ 
    |\gamma-\gamma_q| > R}\hspace{-7mm} \exp\left(-  I(0,\gamma/\eps) -   I^+(0,-\gamma/\eps)\right)  ~+~
  (1+\delta') c.   
\end{align*}
Remark finally that 
\[
 I(0,\gamma/\eps) +  I^+(0,-\gamma/\eps) ~=~ \Bigl(I(0,\gamma/|\gamma|) +
 I^+(0,-\gamma/|\gamma|)\Bigr) |\gamma|/\eps ~\geq~ I_{min} |\gamma|/\eps 
\]
where by Lemma~\ref{lem6-1}, 
\[
I_{min} ~\dot=~ \inf_{\gamma\in\R^{d-1}\times\{0\} :~|\gamma|=1} I(0,\gamma) + I^+(0,-\gamma)
~>~ 0.   
\]
This proves that the right hand side of \eqref{e6-7} does not
exceed 
\begin{align*}
&\hspace{-0.5cm}\limsup_{\eps\to 0} ~\eps
  \log \hspace{-2mm} \sum_{n\geq R - |\gamma_q |} \hspace{-1.5mm} \text{Card}\{\gamma \in
    \eps \Z^d :~ 
    n\leq |\gamma|\leq n+1\} \exp(-  I_{min} n/\eps)  +
  (1+\delta') c \\    
&\leq \limsup_{\eps\to 0} ~\eps
  \log \hspace{-2mm} \sum_{n\geq R - |\gamma_q |} (1 + 2(n + 1)/\eps)^{d-1}  \exp(-  I_{min} n/\eps)  ~+~
  (1+\delta') c \\
&\leq - I_{min} (R-|\gamma_q|) + (1+\delta') c.
\end{align*}
The inequality \eqref{e6-7} holds therefore for $R > |\gamma_q| + ((1+\delta') c  + I(0,q))/I_{min} $.
\end{proof}

\section{Ratio limit theorem for Markov-additive processes}\label{sec-rlth}
In this section we recall the ratio limit theorem for Markov-additive processes. 
 
A Markov chain 
${\cal Z}(t)=(A(t),M(t))$ on $\Z^{d-1}\times \N$ with transition probabilities  
  $p\bigl((x,y),(x',y')\bigr)$ is called {\em Markov-additive} if 
\[
p\bigl((x,y),(x',y')\bigr) ~=~ p\bigl((0,y),(x'-x,y')\bigr)
\]
for all $x,x'\in\Z^{d-1}$, $y,y'\in \N$. $A(t)$ is an {\em additive} part of the process
${\cal Z}(t)$,  and $M(t)$ is its {\em Markovian part}. The Markovian part $M(t)$ is a Markov
chain on $\N$ with transition  probabilities 
\[
p_M(y,y') = \sum_{x\in\Z^{d-1}} p\bigl((0,y),(x,y')\bigr). 
\]

The assumption we need on the Markov-additive process ${\cal Z}(t)=(A(t),M(t))$ are the following. 

\begin{enumerate}
\item[(A1)] There exist 
$\theta >0$ and $C>0$ such that for any $z,z'\in\Z^{d-1}\times \N$ there is a sequence of
  points $z_0, z_1, \ldots,z_n\in\Z^{d-1}\times \N$  with $z_0=z$, $z_n=z'$ and 
    $n\leq C|z'-z|$ such that  
\[
 |z_i-z_{i-1}| \leq C \quad \text{ and } \quad \P_{z_{i-1}}({\cal Z}(1) = z_i) \geq \theta,
 \quad \quad \forall \;
 i=1,\ldots,n. 
\]
\item[(A2)] {\em The function 
\[
\hat\varphi(a) ~=~ \sup_{z\in \Z^{d-1}\times \N} 
~\E_{z}\bigl( \exp(a\cdot ({\cal Z}(1)-z))\bigr)   
\]
is finite everywhere on $\R^{d}$. }
\item[(A3)] {\em Up to multiplication by constants, there is a unique 
 positive harmonic function   $h$ of the Markov process ${\cal Z}(t)=(A(t),M(t))$ such that 
\be\label{e7-1}
\sup_{x\in\Z^{d-1}}h(x,y) < \infty.
\ee}
\end{enumerate} 
Remark that the Markov-additive process ${\cal Z}(t)=(A(t),M(t))$ is not necessarily
stochastic~: in some points $z = (x,y)\in\Z^{d-1}\times \N$, 
the transition matrix can be strictly sub-stochastic. When
the  Markov-additive 
process ${\cal Z}(t)=(A(t),M(t))$ is stochastic,  the last assumption means that  the only
positive harmonic functions  $h : \Z^{d-1}\times\N\to\R_+$ satisfying \eqref{e7-1} 
are constant.

If the assumption (A1) is satisfied then there is a bounded function
$n_0 : \N \to \N^*$ such that for any $z=(x,y)\in\Z^{d-1}\times \N$, 
\[
p^{(n_0(y))}\bigl((x,y), (x,y)\bigr) \geq \theta^{n_0(y)} > 0
\]
and hence, there is $k  \in\N^*$  (for instance, $k = n!$ with $n = \max_y n_0(y)$) such 
that 
\[
 p^{(k)}\bigl(z, z\bigr) \geq  \theta^{k}, \quad \quad \quad \forall z\in\Z^{d-1}\times\N.
\]
Let $\hat{k}$ be the greatest common divisor of the set of all integers $k >0$ for
which 
\[
\inf_{z\in Z^{d-1}\times\N} p^{(k)}(z,z) ~>~ 0
\] 
then from (A3) it follows that 
\begin{enumerate}
\item[(A3')] {\em Up to multiplication by constants, there is a unique 
 positive harmonic function   $h$ of the Markov process ${\cal Z}(t)=(A(t),M(t))$ satisfying the equality  
$h(z + \hat{k}w) = h(z)$ for all $z\in\Z^{d-1}\times\N$ and $w\in
 \Z^{d-1}\times\{0\}$.}
\end{enumerate}

We will use the following property of Markov-additive processes. $G(z,z')$ denotes here
Green's function of the Markov process ${\cal Z}(t)=(A(t),M(t))$.

\begin{prop}\label{pr7-1} Let a Markov-additive process ${\cal
    Z}(t)=(A(t),M(t))$ be 
  transient and satisfy the hypotheses (A1), (A2), (A3).  Suppose moreover that a sequence of
  points $z_n\in\Z^{d-1}\times \N$ 
  is such that  $|z_n|\to\infty$ and  
\[
\liminf_{n\to\infty}~\frac{1}{|z_n|} \log G\bigl(z_0, z_n\bigr)  ~\geq~ 0.
\]
Then 
\[
\lim_{n\to\infty} G(z,z_n)/G(z',z_n) ~=~ h(z)/h(z') 
\]
for all $z,z'\in\Z^{d-1}\times E$.
\end{prop}

For a Markov-additive processes ${\cal Z}(t)=(A(t),M(t))$ with a one-dimensional additive part and for
$z_n=(n,y)$ with a given $y\in \N$, this property was obtained by Foley and
McDonald~\cite{Foley-McDonald}. In the present setting, under the hypotheses (A1), (A2)
and (A3'), the proof of this proposition is given
in~\cite{Ignatiouk:06}. 

\section{Proof of Theorem~\ref{th1-1}}\label{sec-proof}
Under the hypotheses (H1)-(H4), the interior of the set $\hat{D} ~\dot=~ 
\{a\in D :~ \varphi_0(a) \leq 1\}$
is non-empty because $\varphi(0) = \varphi_0(0) = 1$, $
\nabla\varphi(0) ~=~ m ~\not=~ 0$ 
and 
\[
\frac{\nabla\varphi(0)}{|\nabla\varphi(0) |} +
\frac{\nabla\varphi_0(0)}{|\nabla\varphi_0(0)|} ~=~ \frac{m}{|m|} + \frac{m_0}{|m_0|} ~\not=~ 0. 
\]
From this it follows that $
Card\left( (\Theta\times\R)\cap \partial_+ D\right) ~>~ 1$ 
because the orthogonal projection of the set $\{a\in D :~ \varphi_0(a) \leq 1\}$ on the
hyper-plane $\R^{d-1}\times\{0\}$ is homeomorphic to the set $(\Theta\times\R)\cap
\partial_+ D$ where 
\[
\Theta ~\dot=~ \{ \a \in\R^{d-1} :~ \inf_{\beta\in\R} \max\{\varphi(\a,\beta),
\varphi_0(\a,\beta)\} \leq 1\}.
\]
 By Proposition~\ref{pr3-1}, this proves  that there are non-constant
non-negative harmonic functions and consequently, by Theorem~6.2 of \cite{Seneta}, the
Markov process $Z(t)$ is 
transient. The first assertion of Theorem~\ref{th1-1} is therefore proved.

To prove the second assertion we have  to show that 
\be\label{e8-1}
\lim_{n\to\infty} G(z,z_n)/G(z_0,z_n) = h_{\hat{a}(q)}(z)/h_{\hat{a}(q)}(z_0),\quad \quad \quad \forall \;
z\in\Z^{d-1}\times\N.  
\ee
for any non-zero vector
   $q\in{\cal S}_+^d$,  and any sequence of points 
  $z_n\in\Z^{d-1}\times\N$ with $\lim_{n\to\infty} |z_n| = +\infty$ and $\lim_{n\to\infty}
z_n/|z_n|= q$.
The proof of \eqref{e8-1} is different in each of the following cases~: 
\noindent
\begin{itemize}
\item[--] Case 1~: $q\in\R^{d-1}\times\{0\}$, 
\item[--] Case 2~: $q\in\R^{d-1}\times\R_+^*$, and 
    $\varphi_0(\ol{\hat{a}(q)}) < 1$, 
\item[--] Case 3~: $q\in\R^{d-1}\times\R_+^*$, \,
   $\varphi_0(\ol{\hat{a}(q)}) = 1$ and  $\gamma_q\not=0$, 
\item[--] Case 4~:  $q\in\R^{d-1}\times\R_+^*$, \, $\varphi_0(\ol{\hat{a}(q)}) = 1$ and  $\gamma_q=0$, 
\end{itemize}

Recall that $a=\hat{a}(q)$ is the only point of the set $(\Theta\times\R)\cap \partial_+
D$ for which $q\in V(a)$ (see Lemma~\ref{lem2-5}). We denote by $V(a)$ 
the normal cone to the set $(\Theta\times\R)\cap D$ at the point $a$. By
Lemma~\ref{lem2-3}, for every $a\in (\Theta\times\R)\cap \partial_+D$,  
\[
V(a) ~=~ V_D(a) + \left( V_D(\ol{a}) + V_{D_0}(\ol{a})\right)\cap
\left(R^{d-1}\times\{0\}\right) 
\]
where   $\ol{a}$ is the only
point in the boundary $\partial_-D = \{a\in\partial D : \nabla\varphi(a)
\in\R^{d-1}\times\R_-\}$ which has the same orthogonal projection to the hyper-plane as
the point $a$, 
$$
V_D(a) = \{ c \nabla\varphi(a)  \; |\; c \geq 0\}$$ 
is the normal cone to the set $D$ at the point $a$ and 
\[
V_{D\cap D_0}(\ol{a})  = V_D(\ol{a}) + V_{D_0}(\ol{a}) = \{c_1\nabla\varphi(\ol{a}) +
c_2\nabla\varphi_0(\ol{a}) \; | \; c_1,c_2 \geq 0\}
\]
is the normal cone to the set $D\cap D_0$ at the point $\ol{a}$. For
$q\in\R^{d-1}\times\R_+^*$, according to
Corollary~\ref{cor2-4},
\[
\gamma_q \in \left( V_D(\ol{\hat{a}(q)}) + V_{D_0}(\ol{\hat{a}(q)})\right)\cap
\left(R^{d-1}\times\{0\}\right) 
\]
is the only vector at the
hyper-plane $\R^{d-1}\times\{0\}$ for which $q-\gamma_q, \gamma_q \in V(\hat{a}(q))$ and
$q-\gamma_q\in V_D(\hat{a}(q))$. By Lemma~\ref{lem2-5}, for $\gamma_q\not= 0$ we have therefore 
\be\label{e8-2}
\hat{a}(q) = \hat{a}(\gamma_q) = \hat{a}(q-\gamma_q). 
\ee
Recall finally that for every $a\in (\Theta\times\R)\cap
\partial_+D$, 
\be\label{e8-3}
\varphi(\ol{a}) =1 \quad \text{and} \quad \varphi_0(\ol{a}) \leq 1
\ee
because  $\ol{a} \in\partial_-D \subset \partial D$ according to the definition of the
mapping $a\to\ol{a}$, and $\ol{a}\in D_0$ according to the definition of the set
$\Theta$.

{\bf Case 1~: }  To get  \eqref{e8-1} in this case we combine  the ratio limit theorem and
the method of the exponential change of measure~: Proposition~\ref{pr7-1} is applied for
a twisted Markov process $\tilde{Z}(t)$ on $\Z^{d-1}\times\N$ having transition 
probabilities 
\begin{align}
\tilde{p}(z,z') &~=~ \quad \exp\left(a\cdot (z'-z)\right) p(z,z') \nonumber \\ &~=~ \begin{cases}
  \exp\left(a\cdot (z'-z)\right)  \mu_0(z'-z) &\text{if
  $z\in\Z^{d-1}\times\{0\}$} \\
  \exp\left(a\cdot (z'-z)\right)  \mu(z'-z) &\text{if
  $z\in\Z^{d-1}\times\N^*$} \\
\end{cases}
\label{e8-4}
\end{align}
with $a = \ol{\hat{a}(q)}$. 
The infinite matrix $(\tilde{p}(z,z'), \;
z,z'\in\Z^{d-1}\times\N)$ is substochastic because $\varphi(\ol{\hat{a}(q)}) = 1$ and
$\varphi_0(\ol{\hat{a}(q)}) \leq 1$ (see \eqref{e8-3}). Green's function
$\tilde{G}(z,z')$ of the twisted Markov process $\tilde{Z}(t)$ 
  satisfies the equality   
\be\label{e8-5}
\tilde{G}(z,z') ~=~ G(z,z') \exp\left(\ol{\hat{a}(q)}\cdot(z'-z)\right), \quad \quad
\forall z,z'\in\Z^{d-1}\times\N 
\ee  
and hence, using Proposition~\ref{pr5-6} we get 
\begin{align}
\liminf_{n\to\infty} \frac{1}{|z_n|} \log  \tilde{G} (z,z_n) &~=~ \ol{\hat{a}(q)}\cdot q +
\liminf_{n\to\infty} \frac{1}{|z_n|} \log  G(z,z_n) \nonumber\\ &~=~ \left(\ol{\hat{a}(q)} -
\hat{a}(q)\right)\cdot q  ~=~ 0 \label{e8-6}
\end{align} 
where the last relation holds because  $q\in\Z^{d-1}\times\{0\}$ and 
 the orthogonal
projections of the points $\ol{\hat{a}(q)}$ and $\hat{a}(q)$  on the hyper-plane
$\R^{d-1}\times\{0\}$ are 
  identical according to the definition of the mapping $a\to \ol{a}$.  Furthermore, we
  have to check that the  
twisted Markov-additive process $\tilde{Z}(t)$ satisfies the hypotheses
 (A1),(A2) and (A3)  of Section~\ref{sec-rlth}. For this we first notice that the  
 Markov process $Z(t)$ satisfies communication conditions (A1) because of 
 Proposition~\ref{pr4-1} :~ for any
$z,z'\in\Z^{d-1}\times\N^*$ there is a sequence of 
  points $z_0, z_1, \ldots,z_n\in\Z^{d-1}\times\N^*$  with $z_0=z$, $z_n=z'$  and 
    $n\leq C|z'-z|$  such that  and 
\[
 |z_i-z_{i-1}| \leq C \quad \text{ and } \quad \P_{z_{i-1}}(Z(1) = z_i) ~\geq~ \theta, \quad \quad \forall \;
 i=1,\ldots,n. 
\] 
For the twisted Markov process $\tilde{Z}(t)$ we have therefore 
\begin{align*}
 \P_{z_{i-1}}(\tilde{Z}(1) = z_i)  &\geq  \exp\left(- \ol{\hat{a}(q)}\cdot (z_i-z_{i-1})\right)
 \theta ~\geq~ \exp\left(- C\left|\ol{\hat{a}(q)}\right|\right) \theta, 
\end{align*}
for all $i=1,\ldots,n$ and consequently,  $\tilde{Z}(t)$ also satisfies communication
condition~(A1). Next, we remark that by Proposition~\ref{pr3-1}, the constant multiples of
 the function $h_{\hat{a}(q)}$  are the only non-negative harmonic functions of the Markov
 process $Z(t)$ for which  
\[
\sup_{x\in\R^{d-1}} \exp(- \hat\a(q) \cdot x)  h(x,y) ~<~ +\infty, \quad \quad \forall
y\in\N 
\] 
where $\hat{\a}(q)$ denotes the $d-1$ first coordinates
of the point $\hat{a}(q)$. The constant multiples of the function 
\[
\tilde{h}(z) = \exp(-\ol{\hat{a}(q)}\cdot z) h_{\hat{a}(q)}(z) 
\]  
are therefore the only
 non-negative harmonic functions of the twisted Markov 
 process $\tilde{Z}(t)$ for which 
\[
\sup_{x\in\R^{d-1}}  \tilde{h}(x,y) ~<~ +\infty \quad \quad \forall
y\in\N. 
\]
Finally, the function  
\[
\sup_{z\in\Z^{d-1}\times\N} \E_z\bigl(\exp(a\cdot (\tilde{Z}(1)-z))\bigr) ~=~
\max\left\{\varphi\bigl(a +\ol{\hat{a}(q)}\bigr), \varphi_0\bigl(a +\ol{\hat{a}(q)}\bigr)\right\} 
\]
is finite everywhere on $\R^d$ because of the assumption (H4).  
The twisted Markov process $\tilde{Z}(t)$ satisfies therefore  the hypotheses
 (A1),(A2) and (A3) of Section~\ref{sec-rlth}. Using Proposition~\ref{pr7-1} together
with \eqref{e8-6} we get  
\[
\lim_{n\to\infty} \tilde{G}(z,z_n)/\tilde{G}(z_0,z_n) =
 \tilde{h}(z)/\tilde{h}(z_0),\quad \quad \quad \forall \; 
z\in\Z^{d-1}\times\N  
\]
and hence, using again \eqref{e8-5} we obtain \eqref{e8-1}. 

\medskip
\noindent
{\bf Case 2 :~} Suppose now that $q\in\R^{d-1}\times\R_+^*$ and
$\varphi_0\left(\ol{\hat{a}(q)}\right) ~<~ 1$. Here, we
apply Proposition~\ref{pr7-1} for a twisted Markov process $\tilde{Z}(t)$ having  
transition probabilities $ 
\tilde{p}(z,z') ~=~ p(z,z') h_{\hat{a}(q)}(z')/ h_{\hat{a}(q)}(z)$ and  Green's function
\be\label{e8-7}
\tilde{G}(z,z') = G(z,z')  h_{\hat{a}(q)}(z')/
h_{\hat{a}(q)}(z). 
\ee
Such a Markov process is usually called  $h$-transform of the original Markov process $Z(t)$. It is
Markov-additive as well as the  Markov process $Z(t)$ 
because the harmonic function $h_{\hat{a}(q)}$ satisfies the equality 
$h_{\hat{a}(q)}(x,y) = h_{\hat{a}(q)}(0,y) \exp(\hat{\a}(q) \cdot x) $ for all
$(x,y)\in\Z^{d-1}\times\N$. Using quite the same arguments as in the previous case one can
easily show that the new Markov-additive process $\tilde{Z}(t)$ 
satisfies the conditions (A1),(A2) and (A3) of Section~\ref{sec-rlth}. The last condition (A3) is
satisfied here with the constant harmonic function $\tilde{h}(z) \equiv 1$.    Moreover,  from the
explicit representation \eqref{e1-10} of the harmonic function $h_{\hat{a}(q)}$  it follows
that 
\[
\lim_{n\to\infty} \frac{1}{|z_n|} \log h_{\hat{a}(q)}(z_n) ~=~ \hat{a}(q)\cdot q
\]
and hence, by Proposition~\ref{pr5-6},   
\[
\liminf_{n\to\infty} \frac{1}{|z_n|} \log  \tilde{G} (z,z_n) ~=~ \hat{a}(q)\cdot q +
\liminf_{n\to\infty} \frac{1}{|z_n|} \log  G(z,z_n)   ~=~ 0.
\] 
Using Proposition~\ref{pr7-1} we conclude therefore that 
\[
\lim_{n\to\infty} \tilde{G}(z,z_n)/\tilde{G}(z_0,z_n) = 1,\quad \quad \quad \forall \; 
z\in\Z^{d-1}\times\N  
\]
and using next \eqref{e8-7} we get \eqref{e8-1}.

\medskip
\noindent
{\bf Case 3 :~} Suppose now that $q\in\R^{d-1}\times\R_+^*$, \, $\varphi_0\left(\ol{\hat{a}(q)}\right) ~=~
  1$ and $\gamma_q \not=0$. Recall that in this case, 
\be\label{e8-8}
h_{\hat{a}(q)}(z) ~=~ \exp\left(\ol{\hat{a}(q)}\cdot z\right), \quad \quad \quad \forall
\; z\in\Z^{d-1}\times\N.   
\ee
Here, we can not use the above arguments  because \eqref{e8-6} does not hold
and there is no harmonic function satisfying the equality \eqref{e8-7}. To prove \eqref{e8-1} for such a vector $q\in\R^{d-1}\times\R_+^*$ 
we use Proposition~\ref{pr6-1} which proves that for any $\delta >0$ and  $z\in\Z^{d-1}\times\N$, 
\be\label{e8-9}
G(z,z_n) ~\sim~ \Xi^q_{\delta}(z,z_n) \;\text{ as $n\to\infty$} 
\ee
where 
\be\label{e8-10}
\Xi^q_{\delta}(z,z_n) ~\dot=~ \sum_{\substack{w\in
    \Z^{d-1}\times\{0\} :~ |w-\gamma_q |z_n||\leq |z_n|\delta, \\ w'\in \Z^{d-1}\times\N^* 
    }}G(z,w)  \mu_0(w') G_+(w+w',z_n)  
\ee
In Case 1~, we have
already proved that  for all 
$z,z_0\in\Z^{d-1}\times\N$,  
\[
G(z,w)/G(z_0,w) \to h_{\hat{a}(\gamma)}(z)/h_{\hat{a}(\gamma)}(z')
\]
when $|w|\to\infty$ and $w/|w|\to \gamma/|\gamma|\in\R^{d-1}\times\{0\}$. For
any $\sigma >0$ there are therefore 
$n_\sigma >0$ and $\delta >0$ such that  
\[
(1-\sigma) h_{\hat{a}(\gamma_q)}(z)/h_{\hat{a}(\gamma_q)}(z_0) ~\leq~ G(z,w)/G(z_0,w) ~\leq~ (1+\sigma)
h_{\hat{a}(q)}(z)/h_{\hat{a}(q)}(z_0)  
\]
whenever $|w- \gamma_q|z_n|| < \delta |z_n|$ and $n>n_\sigma$. Using these inequalities in
\eqref{e8-10} we obtain 
\[
(1-\sigma) \frac{h_{\hat{a}(\gamma_q)}(z)}{h_{\hat{a}(\gamma_q)}(z_0)} ~\leq~ \frac{ \Xi^q_{\delta}(z,z_n)}{\Xi^q_{\delta}(z_0,z_n)} ~\leq~
  (1+\sigma) \frac{h_{\hat{a}(\gamma_q)}(z)}{h_{\hat{a}(\gamma_q)}(z_0)}. 
\]
for all $n>n_\sigma$. Next, letting $n\to\infty$ and using \eqref{e8-9} we get 
\begin{align*}
(1-\sigma)
  \frac{h_{\hat{a}(\gamma_q)}(z)}{h_{\hat{a}(\gamma_q)}(z_0)} &\leq~  \liminf_{n\to\infty} \frac{
  \Xi^q_\delta(z,z_n)}{\Xi^q_\delta(z_0,z_n)} ~=~ \liminf_{n\to\infty} \frac{G(z,z_n)}{G(z_0,z_n)}   \\
&\leq~ \limsup_{n\to\infty} \frac{G(z,z_n)}{G(z_0,z_n)} ~=~ \limsup_{n\to\infty} \frac{\Xi^q_\delta(z,z_n)}{\Xi^q_\delta(z_0,z_n)} ~\leq~ (1+\sigma)
  \frac{h_{\hat{a}(\gamma_q)}(z)}{h_{\hat{a}(\gamma_q)}(z_0)} 
\end{align*}
and finally, letting $\sigma\to 0$ we conclude that  
\[
\lim_{n\to\infty} {G(z,z_n)}/{G(z_0,z_n)} ~=~ {h_{\hat{a}(\gamma_q)}(z)}/{h_{\hat{a}(\gamma_q)}(z_0)}.
\]
The last relation combined with \eqref{e8-2} proves \eqref{e8-1}. 

\medskip
\noindent
{\bf Case 4 :~ } Suppose finally that $q\in\R^{d-1}\times]0,+\infty[$, \, $\varphi_0\left(\ol{\hat{a}(q)}\right) ~=~
  1$ and
  $\gamma_q =0$.     Here, the harmonic function $h_{\hat{a}(q)}$ is
  defined by \eqref{e8-8}. 
Since in this case 
$\varphi_0\left(\ol{\hat{a}(q)}\right) ~=~ \varphi\left(\ol{\hat{a}(q)}\right) ~=~ 1$ then
  without any
restriction of generality we can assume that 
\be\label{e8-11}
\ol{\hat{a}(q)} = 0. 
\ee 
Otherwise, all the arguments below can be applied for the twisted
Markov process having transition probabilities \eqref{e8-4} with $a= \ol{\hat{a}(q)} $.
So to get \eqref{e8-1} we have to prove  that 
\be\label{e8-12}
\lim_{n\to\infty} G(z,z_n) / G(z',z_n) ~=~ 1 \quad \text{for all
  $z,z'\in\Z^{d-1}\times\N$}. 
\ee
We first prove this relation for the case when $z'-z\in\Z^{d-1}\times\{0\}$. 
For this we combine 
  Proposition~\ref{pr6-1} and the results of Ignatiouk-Robert~\cite{Ignatiouk:06}. Recall
  that 
\[
G(z,z_n) ~=~ G_+(z,z_n) + \hspace{-4mm}\sum_{w\in
    \Z^{d-1}\times\{0\}, w'\in \Z^{d-1}\times\N^*}\hspace{-8mm} G(z,w)  \mu_0(w') G_+(w + w',z_n)
\]
where $G_+(z,z')$ is Green's function of
the homogeneous random walk $Z_+(t)$ on $\Z^{d-1}\times\N^*$ having transition
probabilities $p(z,z') = \mu(z-z')$ which is  killed upon hitting the boundary hyper-plane 
  $\Z^{d-1}\times\{0\}$. 
By Proposition~\ref{pr6-1}, when $n\to\infty$, 
\[
G(z,z_n) \sim 
\Xi_{\delta}^q(z,z_n) ~\dot=~ G_+(z,z_n) + \hspace{-4mm}\sum_{\substack{w\in
    \Z^{d-1}\times\{0\}, w'\in \Z^{d-1}\times\N^* :\\ |w|\leq \delta|z_n|
    }}\hspace{-8mm} G(z,w)  \mu_0(w') G_+(w+w',z_n)
\]
and for any $z'=z+u$ with $u\in\Z^{d-1}\times\{0\}$, 
\[
G(z',z_n) = G(z,z_n-u) \sim 
\Xi_{\delta}^q(z,z_n-u) 
\]
where 
\begin{align}
\Xi_{\delta}^q(z,z_n-u) 
&~\dot=~ G_+(z,z_n-u) + \hspace{-4mm}\sum_{\substack{w\in
    \Z^{d-1}\times\{0\}, w'\in \Z^{d-1}\times\N^* :\\ |w|\leq \delta|z_n|
    }}\hspace{-8mm} G(z,w)  \mu_0(w') G_+(w+w',z_n-u) \nonumber \\
&=~ G_+(z+ u, z_n) + \hspace{-4mm}\sum_{\substack{w\in
    \Z^{d-1}\times\{0\}, w'\in \Z^{d-1}\times\N^* :\\ |w|\leq \delta|z_n|
    }}\hspace{-8mm} G(z,w)  \mu_0(w') G_+(w'+u,z_n-w) \label{e8-13}
\end{align}
Theorem~1 combined with Proposition~2.1 of Ignatiouk~\cite{Ignatiouk:06}
proves that for all $w_0,w''\in\Z^{d-1}\times\N^*$, 
\[
\frac{G_+(w'',v)}{G_+(w_0,v)} \to \frac{\exp(a(q)\cdot w'') - \exp(\ol{a(q)}\cdot w'')}{\exp(a(q)\cdot w_0) -
exp(\ol{a(q)}\cdot w_0)}
\]
as $|v|\to\infty$ and $v/|v|\to q\in\R^{d-1}\times]0,+\infty[$,
    $v\in\Z^{d-1}\times\N^*$. Recall that  $a(q)$ denotes the unique point on the boundary
    $\partial D$ of the set $D=\{a :~\varphi(a)\leq 1\}$ where the vector $q$ is normal to
    $D$.  In our case $q=q -\gamma_q$ and by Corollary~\ref{cor2-4}, the vector $q -\gamma_q$ is normal to
    the set $D$ at the point $\hat{a}(q)$. Hence $a(q) ~=~ \hat{a}(q)$ and according to
    our assumption \eqref{e8-11}, 
\[
\ol{a(q)} ~=~ \ol{\hat{a}(q)} ~=~ 0, 
\]
from which it follows that 
\[
G_+(w'',v)/G_+(w_0,v) \to (\exp(\hat{a}(q)\cdot w') - 1)/(\exp(\hat{a}(q)\cdot
  w_0) - 1)
\]
as $|v|\to\infty$ and $v/|v|\to q\in\R^{d-1}\times]0,+\infty[$,
    $v\in\Z^{d-1}\times\N^*$. In particular, for 
    $u\in\Z^{d-1}\times\{0\}$, from the definition of the mapping $a\to\ol{a}$ it follows
    that  
\[
\hat{a}(q)\cdot u = \ol{\hat{a}(q)}\cdot u = 0
\] 
and consequently, 
\be\label{e8-14}
\lim_{n\to\infty} \frac{G_+(z+ u,z_n)}{G_+(z,z_n)} ~=~ 1, \quad \quad \forall
z\in\Z^{d-1}\times\N^*. 
\ee
Moreover, by Lemma~4.1 of Ignatiouk~\cite{Ignatiouk:06}, the Markov process
$(Z_+(t))$ satisfies communication condition on $\Z^{d-1}\times\N^*$ ~: there exist  
$0< \theta < 1$ and $C>0$ such that for any $w_0,w''\in\Z^{d-1}\times\N^*$ there is a sequence of
  points $w_1, \ldots,w_n\in\Z^{d-1}\times\N^*$  with $w_n=w''$ and 
    $n\leq C|w''-w_0|$ such that  
\[
|w_i-w_{i-1}| \leq C \quad \text{ and } \quad 
\mu(w_i-w_{i-1}) \geq \theta,  \quad \forall \;
 i=1,\ldots,n. 
\]   
The probability that the Markov process $Z_+(t)$
starting at $w_0$ ever hits the point $w''$ is therefore greater than $\theta^{n} \geq
\theta^{C|w_0-w''|}$ which implies that 
\[
G_+(w',v)/G_+(w_0,v) \leq \theta^{-C|w_0-w''|}
\]
 for all $v,w'',w_0\in\Z^{d-1}\times\N^*$. Since the exponential functions are integrable
 with respect to the probability measure $\mu_0$, by dominated convergence theorem from this it
 follows that 
\be\label{e8-15}
\sum_{w'\in\Z^{d-1}\times\N^*} \mu_0(w') \frac{G_+(w'+u,v)}{G_+(w_0,v)} ~\to~ 
 \sum_{w'\in\Z^{d-1}\times\N^*} \mu_0(w') \frac{\exp(\hat{a}(q)\cdot
  (w'+u)) - 1}{\exp(\hat{a}(q)\cdot
  w_0) - 1} 
\ee
as $|v|\to\infty$ and $v/|v|\to q\in\R^{d-1}\times]0,+\infty[$,
    $v\in\Z^{d-1}\times\N^*$. Remark finally that the right hand side of the above display
    is equal to 
\[
 \sum_{w'\in\Z^{d-1}\times\N} \mu_0(w') \frac{\exp(\hat{a}(q)\cdot
  w') - 1}{\exp(\hat{a}(q)\cdot
  w_0) - 1}  ~=~ \frac{\varphi_0(\hat{a}(q)) - 1}{\exp(\hat{a}(q)\cdot
  w_0) - 1}
\]
because $\ol{\hat{a}(q)} =0$ and according to the definition of the mapping $a\to\ol{a}$,
\[\hat{a}(q)\cdot w = \ol{\hat{a}(q)}\cdot w
 \] 
for all $w\in\Z^{d-1}\times\{0\}$.  
Using therefore \eqref{e8-14} and \eqref{e8-15} with $v=z_n-w$ for the right
hand side of \eqref{e8-13} we obtain 
\[
\Xi_{\delta}^q(z, z_n-u) ~\sim~ G_+(z,z_n) + \hspace{-4mm}\sum_{w\in
    \Z^{d-1}\times\{0\} :~ |w|\leq \delta|z_n|
    }\hspace{-8mm} G(z,w) G_+(w_0,z_n-w) \frac{\varphi_0(\hat{a}(q)) - 1}{\exp(\hat{a}(q)\cdot
  w_0) - 1} 
\]
when $n\to\infty$ and $\delta\to 0$. Since the right hand side of the last display does
not depend on $u\in\Z^{d-1}\times\{0\}$ this proves that 
\[
\lim_{\delta\to 0} \liminf_{n\to\infty} \frac{\Xi_{\delta}^q(z, z_n-u)}{
\Xi_{\delta}^q(z, z_n)} ~=~ \lim_{\delta\to 0} \limsup_{n\to\infty} \frac{\Xi_{\delta}^q(z, z_n-u)}{
\Xi_{\delta}^q(z, z_n)} ~=~ 1 
\]
for all
$u\in\Z^{d-1}\times\{0\}$. The equality \eqref{e8-12} for $z\in\Z^{d-1}\times\N$ and
$z'=z+u$ with $u\in\Z^{d-1}\times\{0\}$ follows now from
Proposition~\ref{pr6-1}. 

\medskip

Next, we prove  \eqref{e8-12} for arbitrary $z,z'\in\Z^{d-1}\times\N$. 
Recall that by Proposition~\ref{pr4-1},  the Markov process
$(Z(t))$ satisfies  communication condition on $\Z^{d-1}\times\N$ and consequently, 
there are $0<\delta <1$ and $C>0$ such that for any
$z,z'\in\Z^{d-1}\times\N$, the probability that  the Markov process
$(Z(t))$ starting at $z$ ever hits the point $z'$ is greater than $\theta^{C|z-z'|}$. From this it
  follows that 
\[
\theta^{C|z-z'|} ~\leq~ G(z,z_n)/G(z',z_n) ~\leq~ \theta^{-C|z-z'|} 
\]
for all $z,z'\in\Z^{d-1}\times\N$ and $n\in\N$. Since under the hypotheses (H4), the
exponential functions are integrable with respect to 
the probability measures $\mu$ and $\mu_0$, by dominated convergence theorem we conclude
that for any sub-sequence $n_k$ for which the sequence of functions   
\[
K_n(z) ~=~ G(z,z_{n_k})/G(z_0,z_{n_k})
\] 
converge point-wise, the limit 
\[
K(z) ~\dot=~ \lim_{k\to\infty} K_{n_k}(z) ~\geq~ e^{-\theta |z-z_0|}
\]
is a harmonic function for $(Z(t))$. Remark now that $K(z_0)=1$ and 
\be\label{e8-16}
K(z+u) ~=~ K(z)  \quad \quad  \quad \quad \forall \;
z\in\Z^{d-1}\times\N, \; u\in\Z^{d-1}\times\{0\} 
\ee
because \eqref{e8-12} is already proved for $z'=z+u$ with $u\in\Z^{d-1}\times\{0\}$.
This implies that  $K(z) = 1$ for all $z\in\Z^{d-1}\times\N$ because 
by Proposition~\ref{pr3-1}, the only non-negative harmonic functions  satisfying the
equality \eqref{e8-16} are the constant multiples of the function $h_{\hat{a}(q)}(z)$. 
These arguments
prove that  the sequence of 
functions $K_n$ converge point-wise to the function $K$ because  the function $K$ does not depend
on the sub-sequence $n_k$. The equality \eqref{e8-12} is therefore proved.

\providecommand{\bysame}{\leavevmode\hbox to3em{\hrulefill}\thinspace}
\providecommand{\MR}{\relax\ifhmode\unskip\space\fi MR }
\providecommand{\MRhref}[2]{%
  \href{http://www.ams.org/mathscinet-getitem?mr=#1}{#2}
}
\providecommand{\href}[2]{#2}


\begin{thebibliography}{10}

\bibitem{Billingsley}
Patrick Billingsley, \emph{Convergence of probability measures}, Wiley series
in probability and mathematical statistics, John Wiley \& Sons Ltd, New York,
1968.

\bibitem{D-B}
V.~M. Blinovski{\u\i} and R.~L. Dobrushin, \emph{Process level large deviations
  for a class of piecewise homogeneous random walks}, The Dynkin Festschrift
  Progr. Probab., vol.~34, Birkh\"auser Boston, Boston, MA, 1994, pp.~1--59.

\bibitem{Borovkov:03}
A.A. Borovkov, \emph{On {C}ramer transform, large deviations in boundary-value
  problems and conditional invariance principle}, Sib.Math.J \textbf{3} (1995),
  493--509.

\bibitem{Borovkov:01}
A.A. Borovkov and A.A Mogulskii, \emph{Large deviations for {M}arkov chains in
  the positive quadrant}, Russian Math. Surveys \textbf{56} (2001), no.~5,
  803--916.

\bibitem{D-Z}
Amir Dembo and Ofer Zeitouni, \emph{Large deviations techniques and
  applications}, Springer-Verlag, New York, 1998.

\bibitem{Foley-McDonald}
Robert D.Foley and David R.McDonald, \emph{Bridges and networks: exact
  asymptotics}, Ann. Appl. Probab. \textbf{15} (2005), no.~1B, 542--586.

\bibitem{Doob:03}
J.~L. Doob, J.~L. Snell, and R.~E. Williamson, \emph{Application of boundary
  theory to sums of independent random variables.}, Contributions to
  probability and statistics, Stanford Univ. Press, Stanford, Calif., 1960,
  pp.~182--197.

\bibitem{D-E-W}
P.~Dupuis, R.~Ellis, and A.~Weiss, \emph{Large deviations for {M}arkov
  processes with discontinuous statistics {I} : {G}eneral upper bounds}, Annals
  of Probability \textbf{19} (1991), no.~3, 1280--1297.

\bibitem{D-E-2}
Paul Dupuis and Richard~S. Ellis, \emph{Large deviations for {M}arkov processes
  with discontinuous statistics. {II}. {R}andom walks}, Probability Theory and
  Related Fields \textbf{91} (1992), no.~2, 153--194.

\bibitem{D-E}
\bysame, \emph{The large deviation principle for a general class of queueing
  systems. {I}}, Transactions of the American Mathematical Society \textbf{347}
  (1995), no.~8, 2689--2751.

\bibitem{W-F}
M.~I. Freidlin and A.~D. Wentzell, \emph{Random perturbations of dynamical
  systems}, second ed., Springer-Verlag, New York, 1998, Translated from the
  1979 {R}ussian original by {J}oseph {S}z\"ucs.

\bibitem{Hennequin}
P.L. Hennequin, \emph{Processus de {M}arkoff en cascade}, Ann. Inst. H.
  Poincar\'e \textbf{18} (1963), no.~2, 109--196.

\bibitem{Ignatiouk:02}
Irina Ignatiouk-Robert, \emph{Sample path large deviations and convergence
  parameters}, Annals of Applied Probability \textbf{11} (2001), no.~4,
  1292--1329.

\bibitem{Ignatiouk:04}
\bysame, \emph{Large deviations for processes with discontinuous statistics},
  2005, pp.~1479--1508.

\bibitem{Ignatiouk:06}
\bysame, \emph{Martin boundary of a killed random walk on a half-space}, Journal of
Theoretical Probability, \textbf{21} 2008, no.~ 1, 35-68. 
 

\bibitem{I-M-S}
I.~A. Ignatyuk, V.~A. Malyshev, and V.~V. Shcherbakov, \emph{The influence of
  boundaries in problems on large deviations}, Uspekhi Matematicheskikh Nauk
  \textbf{49} (1994), no.~2(296), 43--102.

\bibitem{Kurkova-Malyshev}
I.A. Kurkova and V.A. Malyshev, \emph{Martin boundary and elliptic curves.},
  Markov Processes Related Fields \textbf{4} (1998), 203--272.

\bibitem{Ney-Spitzer}
P.~Ney and Spitzer F, \emph{The martin boundary for random walk}, Trans. Amer.
  Math. Soc. (1966), no.~121, 116--132.

\bibitem{R}
R.~Tyrrell Rockafellar, \emph{Convex analysis}, Princeton University Press,
  Princeton, NJ, 1997, Reprint of the 1970 original, Princeton Paperbacks.

\bibitem{PhilippeRobert}
Philippe Robert, \emph{Stochastic Metworks and Queues}, Springer-Verlag,
  Berlin, 2003.

\bibitem{Seneta}
E.~Seneta, \emph{Nonnegative matrices and {M}arkov chains}, second ed.,
  Springer-Verlag, New York, 1981.

\bibitem{S-W}
A.~Shwartz and A.~Weiss, \emph{Large deviations for performance analysis},
  Stochastic Modeling Series, Chapman \& Hall, London, London, 1995.

\bibitem{Woess}
Wolfgang Woess, \emph{Random walks on infinite graphs and groups}, Cambridge
  University Press, Cambridge, 2000.

\end{thebibliography}
\end{document}